\documentclass[preprint]{elsarticle}

\usepackage{latexsym}
\usepackage{amsmath}
\usepackage{amssymb}
\usepackage{amsthm}
\usepackage{mathrsfs}
\usepackage{graphicx}
\usepackage{pgfpages}
\usepackage{ifthen}
\usepackage{leftidx,tensor}
\usepackage[T1]{fontenc}
\usepackage[latin1]{inputenc}
\usepackage{mathtools}
\usepackage{comment}
\usepackage{dsfont}

\usepackage[shortlabels]{enumitem}
\usepackage{aliascnt}
\usepackage[bookmarks=true,pdfstartview=FitH, pdfborder={0 0 0}, colorlinks=true,citecolor=red, linkcolor=blue]{hyperref}
\usepackage{bbm}
\usepackage{nicefrac}
\usepackage{romanbar}
\newtheorem{thm}{Theorem}[section]
\newaliascnt{cor}{thm}
\newaliascnt{prop}{thm}
\newaliascnt{lem}{thm}
\newtheorem{cor}[cor]{Corollary}
\newtheorem{prop}[prop]{Proposition}
\newtheorem{lem}[lem]{Lemma}
\aliascntresetthe{cor}
\aliascntresetthe{prop}
\aliascntresetthe{lem} 
\newaliascnt{defn}{thm}
\newaliascnt{asu}{thm}
\newaliascnt{con}{thm}
\aliascntresetthe{defn}
\aliascntresetthe{asu}
\aliascntresetthe{con}
\newcounter{stp}
\newcounter{stpi}
\newcounter{stpci}
\newcounter{stpiii}

\newaliascnt{rem}{thm}
\newaliascnt{exa}{thm}
\newaliascnt{masu}{thm}
\newaliascnt{nota}{thm}
\newaliascnt{sett}{thm}
\newtheorem{rem}[rem]{Remark}

\aliascntresetthe{rem}
\aliascntresetthe{exa}
\aliascntresetthe{masu}
\aliascntresetthe{nota}
\aliascntresetthe{sett}
\numberwithin{equation}{section}

\setlist[enumerate]{font = \normalfont}

\newcommand {\N}	{\mathbb{N}}

\newcommand {\R}	{\mathbb{R}}

\newcommand {\E}	{\mathbb{E}}

\newcommand {\T}	{\mathbb{T}}

\renewcommand{\d}{\, \mathrm{d}}

\DeclareMathOperator{\divH}{div_{\H}}

\newcommand{\sP}{\mathcal{P}}

\renewcommand{\Re}{\mathrm{Re}}

\renewcommand{\H}{\mathrm{H}}

\newcommand{\sigmabar}{\bar{\sigma}}

\newcommand{\ocn}{\mathrm{o}}

	\newcommand{\Gau}{\Gamma_u^{\ocn}}
	\newcommand{\Gab}{\Gamma_b^{\ocn}}

	\newcommand{\dk}[1]{\partial_{#1}}
	\newcommand{\dt}{\dk{t}} 
	\newcommand{\dz}{\dk{z}} 
	
	\renewcommand{\phi}{\varphi}
	
	\renewcommand{\bar}[1]{\overline{#1}}
	\newcommand{\vbar}{\bar{v}}
	
	\newcommand{\vtilde}{\tilde{v}}

	\newcommand{\nablaH}{\nabla_{\H}}
	\newcommand{\DeltaH}{\Delta_{\H}}

	\newcommand{\rC}{\mathrm{C}}
	\newcommand{\rL}{\mathrm{L}}
	
	\newcommand{\rH}{\H}
	\newcommand{\rB}{\mathrm{B}}

	\newcommand{\rE}{\mathrm{E}}

	\newcommand{\rLp}{\rL^p}

	\renewcommand{\P}{\mathcal{P}}

	\renewcommand{\div}{\mathrm{div} \,}
	
	\newcommand{\I}{\, \Romanbar{1} \, }
	\newcommand{\II}{\, \Romanbar{2}\, }
	\newcommand{\III}{\, \Romanbar{3}\, }
	\newcommand{\IV}{\, \Romanbar{4}\, }
	\newcommand{\V}{\, \Romanbar{5}\, }
	\newcommand{\VI}{\, \Romanbar{6}\, }
	\newcommand{\VII}{\, \Romanbar{7}\, }
	\newcommand{\VIII}{\, \Romanbar{8}\, }

\newenvironment{resume}[1]
{\par\medskip
	\noindent\textbf{#1.}\ \itshape}
{\par\medskip}

	\title{The Ekman spiral in primitive equations with wind-driven boundary conditions and Coriolis force}
\hypersetup{pdfauthor={Tim Binz}}

\author[1]{Tim Binz
	}
\ead{tim.binz@princeton.edu}
\affiliation[1]{organization={Program in Applied and Computational Mathematics,
		Princeton University},
	addressline={Fine Hall, Washington Road},
	city={Princeton},
	postcode={08544},
	state={New Jersey},
	country={USA}}
	
\begin{document}


\begin{abstract}
	In this article we investigate the relationship 
	of the 3D incompressible, primitive equations 
	with wind-driven boundary conditions and Coriolis force and the Ekman spiral. We identify three distinct regimes. In the first one, we show that every solution converges exponentially fast to the Ekman spiral as $t \to + \infty$. In particular, this implies that the Ekman spiral is the unique equilibrium of the system.
	In the second regime, we prove nonlinear stability of the Ekman spiral, and in the third one, we show nonlinear instability of the Ekman spiral.  
\end{abstract}

\maketitle	

\begin{resume}{Résumé}
	Dans cet article, nous \'etudions la relation entre les \'equations primitives incompressibles 3D avec des conditions de bord induites par le vent, et avec la force de Coriolis et la spirale d'Ekman.
	Nous identifions trois r\'egimes distincts.
	Dans le premier, nous montrons que chaque solution converge de mani\`ere exponentiellement rapide vers la spirale d'Ekman lorsque $t \to +\infty$.
	En particulier, cela implique que la spirale d'Ekman est l'unique \'equilibre du syst\`eme.
	Dans le deuxi\`eme r\'egime, nous prouvons la stabilit\'e de la spirale d'Ekman, et dans le troisi\`eme, nous montrons l'instabilit\'e de la spirale d'Ekman. 
\end{resume}

\section{Introduction}
\label{sec:intro_2}

In this paper we study a classic scenario: the effect of wind blowing over the ocean (at large scales). 
The primitive equations are the standard model for oceanic
dynamics in large scales \cite{Ped:87_2,LTW:92b_2}. They are derived from the Navier-Stokes equations by assuming the hydrostatic balance for the pressure term.
The 3D incompressible \emph{primitive equations} with gravity in a rotational frame are given by
\begin{equation}
	\left\{
	\begin{aligned}
		\partial_t v - \nu_{\rH} \DeltaH v - \nu_z \partial_z^2 v + u \cdot \nabla v + \frac{1}{\rho_0}\nablaH \pi + f v^\perp &= 0, &&\text{ in } (0,T) \times \Omega, \\
		\partial_z \pi + \rho_0 g &= 0, &&\text{ in } (0,T) \times \Omega, \\
		\div u &= 0, &&\text{ in } (0,T) \times \Omega ,
	\end{aligned}
	\right. 
	\tag{PE}
	\label{eq:pe_2}
\end{equation}
on a periodic layer $\Omega := \T^2 \times (-h,0)$ with boundary $\Gamma = \Gamma_u \cup \Gamma_b$, where $\Gamma_u := \T^2 \times \{ 0 \}$ and $\Gamma_b := \T^2 \times \{ -h \}$.
Here $u = (v,w)$ denotes the velocity of the fluid, where $v = (v_1,v_2)$ is the horizontal component of the velocity and $w$ is the vertical one. Further, $\pi$ is the pressure of the fluid and $\nu_{\rH} > 0$ and $\nu_z > 0$ are the horizontal and vertical viscosities of the fluid. 
Moreover, $f \in \R$ is the $f$-plane approximation of the Coriolis force, $v^\perp := (-v_2,v_1)^\top$, $\rho_0>0$ the reference density, $g > 0$ the gravity constant and $h>0$ the depth.
Finally, $\nablaH g= (\partial_x g,\partial_y g)^\top$ denotes the horizontal gradient, $\divH G = \partial_x G_x + \partial_y G_y$ the horizontal divergence and $\DeltaH := \divH \nablaH$ the horizontal Laplacian. 
The \emph{wind-driven boundary conditions} at the upper surface describes the effect of wind blowing over the surface of the ocean
\begin{equation}
	\dz v = \tau, \quad \text{ on } (0,T) \times \Gamma_u, 
	\label{eq:wind bc_2}
\end{equation} 
where $\tau \in \R^2$ is the given (constant) wind stress.
The lower boundary conditions connect our problem to the region of deep water, where the geostrophic balance holds. We call it the \emph{geostrophic boundary conditions} and it reads as 
\begin{equation}
	v = v_g , \quad \text{ on } (0,T) \times \Gamma_b, 
	\label{eq:bottom_2}
\end{equation}
where $v_g \in \R^2$ is the (constant) geostrophic flow. Further, we assume that there is no vertical velocity at the boundaries, i.e.
\begin{equation}
	w = 0, \quad \text{ on } (0,T) 
	\times \Gamma . 
	\label{eq:w bc_2}
\end{equation}
Finally, this system is completed by the initial data $v(0) = v_0$. 
The 
\emph{Ekman spiral} occurs when wind drags on the surface of the ocean. The Ekman spiral
$v_{\rE}$ is the unique solution of the balance between the Coriolis force and the vertical dissipation, i.e. $v_{\rE}(x,y,z) := v_{\rE}(z)$ solves
\begin{equation}
	\nu_z \cdot \partial_z^2 v_{\rE} = f v_{\rE}^\perp \quad \text{ in } \Omega ,
	\label{eq:ekman_2}
	\tag{Ekman} 
\end{equation}
supplemented by the wind-driven boundary conditions \eqref{eq:wind bc_2} and geostrophic boundary conditions \eqref{eq:bottom_2}.


\medskip 

The Ekman spiral was discovered by Ekman \cite{Ekm:05_2} in 1905. The question of asymptotic stability of the Ekman spiral in the context of Navier-Stokes equations has become an active area of research: 
In \cite{GM:97_2} Grenier and Masmoudi studied the Ekman layer for well prepared initial data.   
This was generalized by Masmoudi in \cite{Mas:00_2} to general data. 
Rousset \cite{Rou:04a_2} established stability of the Ekman layer for sufficiently large rotation in the regime of high regularity. 
In \cite{HHMS:10_2} strong asymptotic stability of the Ekman spiral has been shown. 
Further, Giga and Saal \cite{GS:15_2} proved uniform exponential stability of the Ekman spiral. Their strategy is based on \cite{GS:13_2}.  
In \cite{GIMMS:07_2} local existence and uniqueness of mild solutions of the Ekman-Navier Stokes equations for non-decaying initial data has been shown. 
Fischer and Saal \cite{FS:13_2} established instability of the Ekman spiral when the Reynolds number satisfies $\Re > 55$, under a spectral assumption supported by numerical evidence from Lilly \cite{Lil:66_2}.
The case of two dimensional Navier-Stokes equations with rotation has been studied in \cite{CDGG:02_2}. 
Moreover, Hieber and Stannat \cite{HS:13_2} proved asymptotic stability in the presence of a stochastic force.  

In \cite{Kob:14_2} weak and strong asymptotic stability of the Ekman layer in the Navier-Stokes-Boussinesq system has been studied. 
In \cite{Rou:05_2} Rousset showed the stability of the Ekman-Hartmann boundary layer in MHD equations.
Moreover, in \cite{DDG:99_2, DG:03_2, Rou:04b_2, DGV:17_2} stability of other boundary layers have been studied. 
For a survey of results and more information about the stability of Ekman spirals we refer to the monographs \cite{CDGG:06_2} and \cite{Kob:14_2}.

\smallskip 

In 1922 the primitive equations were introduced by Richardson \cite{Ric:22_2}. 
Their analysis started with a series of papers of Lions, Temam and Wang \cite{LTW:92a_2,LTW:92b_2,LTW:93_2,LTW:95_2}. 
The primitive equations can be derived from the Navier-Stokes equations as a small ratio limit as rigorously shown by Li and Titi \cite{LT:19_2, LTG:22_2} and by Furukawa et al. \cite{FGHHKW:20_2,FGK:21_2,FGHHKW:25_2}. 
The breakthrough result of Cao and Titi \cite{CT:07_2} established global well-posedness of the primitive equations with homogeneous Neumann boundary conditions for large $\rH^1$-initial data. Since then different proofs, refinements and extensions of this result have been found, see e.g.~\cite{Kob:07_2,KZ:07_2,HK:16_2,GGHHK:20_2,GHK:17_2, HHK:16_2, CLT:15_2, CLT:14a_2, CLT:14b_2, CLT:17_2, BH:22_2}.
The role of the rotation in the context of primitive equations has been studied in \cite{GILT:22_2}. 

\smallskip 

Wind-driven boundary conditions for the Navier-Stokes equations have been investigated e.g. in
\cite{DG:00_2,BS:01_2,BGMR:03_2,DS:09_2}. 
The study of wind-driven boundary conditions for the primitive equations goes back to the original articles of Lions, Temam and Wang \cite{LTW:93_2,LTW:95_2}, where they introduced the coupled atmosphere-ocean system consisting of two primitive equations coupled by non-linear wind-driven boundary conditions and showed existence of weak solutions. Global strong well-posedness of the coupled atmosphere-ocean system has been established recently in \cite{BBHZ:25_2}. 
Moreover, primitive equations with stochastic wind-driven boundary conditions have been studied in \cite{BHHS:24_2}. 
For general background about primitive equations and Ekman spirals we refer to the book of Pedlosky \cite{Ped:87_2}.
%

\medskip

Taking into account that the Ekman spiral occurs classically in oceanic dynamics (and sometimes also in atmospheric dynamics) and the primitive equations serve as a standard model for oceanic (and atmospheric)
dynamics it is quite surprising that the stability of the Ekman spiral in the primitive equations has not been studied so far. 

\medskip 

In the present article we investigate the role of the Ekman spiral in the primitive equations \eqref{eq:pe_2} subject to the wind-driven conditions \eqref{eq:wind bc_2} and the geostrophic boundary conditions \eqref{eq:bottom_2}.
The behaviour of the Ekman spiral depends highly on the given parameters: the viscosities/ Reynolds numbers, the strength of the Coriolis force, the depth of the ocean, as well as the wind stress $\tau$ and the geostrophic flow $v_g$. 
In our main theorem, \autoref{thm:main_2}, we identify three distinct parameter regimes exhibiting qualitatively different dynamical behaviour. These regimes are determined by operator-theoretic constants for the linearization of the primitive equations around the Ekman spiral: the \emph{hydrostatic Ekman-Stokes operator} $A$. We denote its spectral bound by $s(A)$, its numerical range by $W(A)$ and define $\alpha := \sup_{\mu \in W(A)} \Re (\mu)$. Note that $s(A) \leq \alpha$ but it is unclear whether equality holds because $A$ is not self-adjoint due to the lower order terms involving the Ekman spiral. 

\begin{itemize}
	\item If $\alpha < 0$ 
	the primitive equations \eqref{eq:pe_2} with \eqref{eq:wind bc_2} and \eqref{eq:bottom_2}
	generate a semi-flow which converges exponentially fast to its unique equilibrium: the Ekman spiral;
	\item If $s(A) < 0$ and $\alpha \geq 0$, the Ekman spiral is stable;
	\item If $s(A) > 0$, the Ekman spiral is unstable. 
\end{itemize}
%
%
%
Loosely speaking this shows that the Ekman spiral in primitive equations becomes more unstable if the Reynolds number increases. The first regime is the laminar regime, the second one the transitional regime and the third one the turbulent regime. 
%
%
%
Our result for the second regime is of a similar flavour to that of Giga and Saal \cite{GS:15_2} about the Ekman spiral in Navier-Stokes equations, whereas the instability result for the third regime is in the spirit of Fischer and Saal \cite{FS:13_2}.
In contrast to these results our result in the first case is non-perturbative: it does not require any smallness conditions of the initial data. Indeed, solutions of \eqref{eq:pe_2} converge to the Ekman spiral for all (!) initial data. This implies, in particular, that the Ekman spiral is the unique equilibrium of the primitive equations, see \autoref{cor:unique_2}.  

\smallskip 

Our proof relies on the analysis of the  \emph{Ekman-primitive equations} \eqref{eq:difference_2} satisfied by the difference $v_d := v - v_{\rE}$ of the solution~$v$ to the primitive equations subject to the wind-driven boundary conditions \eqref{eq:wind bc_2} and the geostrophic boundary conditions \eqref{eq:bottom_2} and the Ekman spiral $v_{\rE}$. 
We apply the maximal regularity theory of Prüss-Wilke \cite{PW:17_2,PWS:18_2} and show $\rL^\infty_t \rH^1$-$\rL^2_t \rH^{2}$-a priori estimates. The main ingredients for this are the maximal regularity of the \emph{hydrostatic Ekman-Stokes operator}, see \autoref{prop:mr_2}, which generalizes the main result from \cite{GGHHK:17_2}, the bilinear estimates from \autoref{lem:bilinear bounds_2}, which are an improvement of the bounds from \cite{GGHHK:20_2}, as well as, the aforementioned a-priori estimates. The location of the spectrum of the hydrostatic Ekman-Stokes operator is sensitive to the lower order terms involving the Ekman spiral. This yields the different results in different regimes of parameters. In the first and the second one the hydrostatic Ekman-Stokes operator has a negative spectral bound whereas in the third the spectral bound is positive. Using this fact part 2 and 3 of the main result follow from both sides of the principle of linearized stability \cite[Section 5]{PS:16_2} or \cite{PSZ:09_2}. 
For the a-priori estimates we follow the strategy of Hieber-Kashiwabara \cite{HK:16_2} using the splitting in baroclinic and barotropic modes and interpolation techniques. In each step here we have to bound additional terms coming from the Ekman spiral. 
Let us comment on these estimates in more detail. 
In the energy estimate, \autoref{lem:energy_2}, we need an additional assumption to the parameter which distinguishes case 1 from 2, to control the terms coming from the Ekman spiral. 
In the  $\rL^2$-estimate of $\partial_z v_d$, \autoref{lem:estimate vz_2}, and the $\rH^1$-estimate of the barotropic mode, \autoref{lem:estimate vbar_2}, these terms can be controlled straightforwardly. In contrast the $\rL^4$-estimates of the baroclinic modes, \autoref{lem:estimate vtilde_2}, require a more delicate analysis. The latter one relies on interpolation techniques and anisotropic Sobolev embeddings.


Additionally, in \autoref{prop:Hk_2}, we establish higher order $\rL^\infty_t \rH^k$-$\rL^2_t \rH^{k+1}$-bounds  for all $k \in \N$ for the difference between the solution to the primitive equations and the Ekman spiral.
For $k = 1$ these bounds are analogous to the classic bounds of Cao and Titi from \cite{CT:07_2}, which leads to the global existence of strong solutions. The $\rL^\infty_t \rH^2$-bound for primitive equations (with homogeneous boundary conditions) also appears the first time in \cite{CT:07_2}, whereas the $\rL^2_t \rH^{3}$-bound has been established recently in \cite[Appendix A]{Fur:24_2}.
To our best knowledge for $k > 2$ these bounds are novel even for primitive equations with homogeneous boundary conditions.
Our proof relies on the maximal regularity of the hydrostatic Ekman-Stokes operator in $\rH^k$-spaces, see \autoref{prop:mr_2}, the improved bilinear estimates from \autoref{lem:bilinear bounds_2} and Gronwall's inequality.
These bounds as well as the method to prove them are of independent interest. 

\medskip

\enlargethispage{\baselineskip}
The present article is organized as follows: 
In \autoref{sec:main_2} we state our main result, \autoref{thm:main_2}, and provide some comments about it. In \autoref{sec:Ekman_2} we study the Ekman spiral of finite depth. 
In \autoref{sec:local_2} we introduce the functional analytic set-up to study our problem. Furthermore, we derive the Ekman-primitive equations and establish the necessary linear theory for the hydrostatic Ekman-Stokes operator in \autoref{prop:mr_2}. Moreover, we provide the novel bilinear estimates of the primitive equations in \autoref{lem:bilinear bounds_2}. 
Combining these facts we conclude the local well-posedness of the Ekman-primitive equations in \autoref{prop:local_2}, and finally the asymptotic stability, \autoref{cor:small data_2}, and instability, \autoref{prop:instability_2}, of the Ekman spiral.
\autoref{sec:global_2} is the central section of this paper. First we establish the blow-up scenario in \autoref{cor:blow-up_2}, which we rule out later by the  $\rL^\infty_t\rH^1 \cap \rL^2_t \rH^2$-estimates of \autoref{prop:estimate_2}. This leads to almost global existence of the Ekman-primitive equations, see \autoref{prop:almost global_2}.
Combining this result with our asymptotic stability of the Ekman spiral from \autoref{cor:small data_2} we conclude case 1 of our main result, \autoref{thm:main_2}: the convergence of every solution of the primitive equations to the Ekman spiral as $t \to +\infty$. 
Moreover, we prove higher-order estimates of the difference in \autoref{prop:estimate_2 Hk} iteratively using our improved bilinear estimates from \autoref{lem:bilinear bounds_2} and the linear theory from \autoref{prop:mr_2}.
Finally, for convenience of the reader we calculate explicitly the Ekman spiral of finite depth in \autoref{sec:ekman explicit_2}.  


\section{Main result}
\label{sec:main_2}

Before stating our main theorem, we introduce the functional-analytic set-up of the primitive equations \eqref{eq:pe_2} with \eqref{eq:wind bc_2} and \eqref{eq:bottom_2}.
Following \cite{CT:07_2}, 
we introduce the \emph{barotropic mode} as the vertical average
\begin{equation*}
	\vbar(t,x,y) := \frac{1}{h} \int_{-h}^{0} v(t,x,y,z) \d z 
\end{equation*}
and define the \emph{baroclinic mode} by $\vtilde := v - \vbar$.
Using $\div u = 0$ and $w = 0$ we obtain
\begin{equation*}
	\divH \vbar = 0  .
\end{equation*}
This observation motivates the definition of the \emph{space of hydrostatic solenoidal $2$-integrable vector-fields} given by 
\begin{equation*}
	\rL^2_{\sigmabar}(\Omega)
	:= \overline{\{ v \in \rC^\infty(\Omega) \colon \divH \vbar = 0 \}}^{\| \cdot \|_{\rL^2}},
\end{equation*}
see \cite{HK:16_2}. 
It yields a Hilbert space when equipped by the canonical scalar product
\begin{equation*}
	\langle v, v' \rangle 
	:= \int_{\Omega} v(x,y,z) \cdot v'(x,x,z) \, \mathrm{d}^3(x,y,z) .
\end{equation*}
Further, we obtain the \emph{hydrostatic Helmholtz projection} $\P \colon \rL^2(\Omega) \to \rL^2_{\sigmabar}(\Omega)$. 
\enlargethispage{\baselineskip}
The key object in our analysis is the linearization of the primitive equations \eqref{eq:pe_2} with boundary conditions \eqref{eq:wind bc_2}--\eqref{eq:w bc_2} around the Ekman spiral \eqref{eq:ekman_2}: 
\emph{the hydrostatic Ekman-Stokes operator} 
$A \colon X_1 \subset X_0 \to X_0$ given by
\begin{equation*}
	\begin{aligned} 
		A v &:= \P ( \nu_{\rH} \DeltaH v +\nu_z \partial_z^2 v)
		- \P (v_{\rE} \cdot \nablaH v 
		+ w \dz v_{\rE}) 
		- \P f v^\perp,
		\\
		X_1 &:= \{ v \in \rH^{2}(\Omega) \cap \rL^2_{\sigmabar}(\Omega) \colon \partial_z v|_{\Gamma_u} = 0, \quad v|_{\Gamma_b} = 0 \} .
	\end{aligned} 
\end{equation*}
We denote the \emph{spectral bound} of $A$ by
\begin{equation*}
	s(A) := \sup_{\lambda \in \sigma(A)} \Re (\lambda) .
\end{equation*}
Further, we denote the \emph{numerical range} of $A$ by
\begin{equation*}
	W(A) := \{ \langle A v , v \rangle \colon v \in D(A) \text{ such that } \| v \|_{\rL^2(\Omega)} = 1 \},
\end{equation*}
and define by
\begin{equation*}
	\alpha := \sup_{\mu \in W(A)} \Re (\mu) . 
\end{equation*} 
Note that $s(A) \leq \alpha$ but equality does not generally hold because $A$ is not self-adjoint due to the lower order terms involving the Ekman spiral.
Both $s(A)$ and $\alpha$ depend on the structure of the Ekman spiral, and thus on the physical parameters of the system: the viscosities $\nu_{\rH}$ and $\nu_z$ (equivalently, the Reynolds numbers), the Coriolis parameter $f$, the depth $h$, as well as the wind stress $\tau$ and the geostrophic flow $v_g$. This dependence is highly nontrivial.

\smallskip   

With this preparation, we are now ready to state our main result.  

\begin{thm}[Dynamical behaviour of the Ekman spiral]\label{thm:main_2}
	Depending on the parameter regime we obtain the following results:
	\begin{enumerate}[(a)]
		\item	
		If $\alpha < 0$ is satisfied,
		then the 3D incompressible primitive equations with Coriolis force \eqref{eq:pe_2} subject to the wind-driven boundary conditions \eqref{eq:w bc_2} and the geostrophic boundary conditions \eqref{eq:bottom_2} are globally wellposed, i.e. for every $v_0 \in \rH^{1}(\Omega) \cap \rL^2_{\sigmabar}(\Omega)$ satisfying \eqref{eq:bottom_2} there exists a unique strong solution $v \in 
		\mathrm{BUC}(\R_+;\rH^{1}(\Omega)\cap \rL^2_{\sigmabar}(\Omega))$ with initial data $v(0) = v_0$ such that the difference $v-v_{\rE} \in\rH^{1}(\R_+;\rL^2(\Omega)) \cap \rL^{2}(\R_+;\rH^{2}(\Omega))$, where $v_{\rE}$ denotes the Ekman spiral \eqref{eq:ekman_2}. 
		The solution $v$ converges exponentially to the Ekman spiral~$v_{\rE}$, so
		\begin{equation*}
			\| v(t) - v_{\rE} \|_{\rH^{1}(\Omega)} \leq C \mathrm{e}^{\omega_0 t}\to 0  \qquad \text{ as } \enspace t \to + \infty , 
		\end{equation*}
		with rate $\omega_0 < 0$.	
		\item If $\alpha \geq 0$ and $s(A) < 0$, then the Ekman spiral $v_{\rE}$ is stable, i.e. there exists $r > 0$ such that for every $v_0 \in \rH^{1}(\Omega) \cap \rL^2_{\sigmabar}(\Omega)$ satisfying \eqref{eq:bottom_2}
		with $\| v_0-v_{\rE} \|_{\rH^1(\Omega)} < r$ there exists 
		a unique global strong solution $v \in \mathrm{BUC}(\R_+;\rH^{1}(\Omega)\cap \rL^2_{\sigmabar}(\Omega))$ with initial data $v(0) = v_0$ such that the difference $v-v_{\rE} \in\rH^{1}(\R_+;\rL^2(\Omega)) \cap \rL^{2}(\R_+;\rH^{2}(\Omega))$ 
		and the solution $v$ converges exponentially to the Ekman spiral~$v_{\rE}$, so
		\begin{equation*}
			\| v(t) - v_{\rE} \|_{\rH^{1}(\Omega)} \leq C \mathrm{e}^{\omega_0 t}\to 0  \qquad \text{ as } \enspace t \to + \infty , 
		\end{equation*}
		with rate $\omega_0 < 0$.
		\item If $s(A) > 0$, then the Ekman spiral $v_{\rE}$ is unstable, i.e. there exists $\varepsilon >0$ such that for any $\delta > 0$ there exists $v_0 \in \rH^{1}(\Omega) \cap \rL^2_{\sigmabar}(\Omega)$ satisfying \eqref{eq:bottom_2} with $\| v_0 - v_{\rE} \|_{\rH^1(\Omega)} < \delta$ such that there is a finite time $t_0 > 0$ with
		\begin{equation*}
			\| v(t_0) - v_{\rE} \|_{\rH^1(\Omega)} \geq \varepsilon . 
		\end{equation*}
	\end{enumerate}	
	%
\end{thm}

Let us now make a few comments on our main result. 

\begin{rem}[Pressure and vertical velocity]\label{rem:w_2}
	The vertical velocity $w$ is determined uniquely by the horizontal velocity $v$ and the geostrophic flow $v_g$ by
	\begin{equation*}
		w(x,y,z)
		= v_g - \int_{-h}^z \divH v(x,y,\xi) \, \d \xi . 
	\end{equation*}
	Moreover, the pressure $\pi$ is determined uniquely (up to constants) by the horizontal velocity $v$ via the elliptic equation
	\begin{equation*}
		\DeltaH \pi 
		= \rho_0 \cdot \biggl(
		-\frac{\nu_z}{h} \divH \bigl(\dz v\big|_{\Gamma_b}\bigr)
		+ \frac{1}{h} \divH \biggl(
		\int_{-h}^0 v \cdot \nablaH v - v \cdot \divH v
		\biggr)
		\biggr) - \rho_0 g \cdot z .
	\end{equation*}
\end{rem}

\begin{rem}[Convergence rate]
	The convergence rate $\omega_0<0$ in \autoref{thm:main_2} is the spectral bound of the linearization, i.e. of the hydrostatic Ekman-Stokes operator $A$ 
	Since the hydrostatic Stokes operator $A$ has compact resolvent $\omega_0$ coincides with the real part of largest eigenvalue.  
\end{rem}

The next result gives precise control over the derivatives of the difference between the solution $v$ to the primitive equations \eqref{eq:pe_2} and the Ekman spiral $v_{\rE}$. 

\begin{prop}[$\rH^k$-bounds]\label{prop:Hk_2}
	Fix $k \in \N$ with $k \geq 2$. Assume $v_0 \in \rH^k(\Omega) \cap \rL^2_{\bar{\sigma}}(\Omega)$ and the boundary conditions \eqref{eq:wind bc_2} and \eqref{eq:bottom_2}.
	Then the $\rH^k$-norms of the difference are controlled by bounds~$B_k$,
	\begin{equation*}
		\| v(t) - v_{\rE} \|_{\rH^{k}(\Omega)}^2
		+ \int_0^t \| v(t) - v_{\rE} \|_{\rH^{k+1}(\Omega)}^2
		\leq B_k(t) .
	\end{equation*} 
\end{prop}

Since \eqref{eq:ekman_2} is merely a linear ordinary differential equation it can be solved explicitly.

\begin{rem}[Ekman Spiral]\label{rem:ekman_2}
	The Ekman spiral $v_{\rE} = (v_1,v_2)^\top$ is explicitly given as
	\begin{align*}
		v_1(z)
		&= k_1 \sin(\tfrac{z}{d}) e^{-\tfrac{z}{d}}
		+ k_2 \cos(\tfrac{z}{d}) e^{-\tfrac{z}{d}}
		+ k_3 \sin(\tfrac{z}{d}) e^{\tfrac{z}{d}}
		+ k_4 \cos(\tfrac{z}{d}) e^{\tfrac{z}{d}}, \\
		v_2(z)
		&=
		k_1 \cos (\tfrac{z}{d}) e^{-\tfrac{z}{d}}
		- k_2 \sin (\tfrac{z}{d}) e^{-\tfrac{z}{d}}
		- k_3 \cos (\tfrac{z}{d}) e^{\tfrac{z}{d}}
		+ k_4 \sin (\tfrac{z}{d}) e^{\tfrac{z}{d}},
	\end{align*}
	where $d := (\frac{2 \nu_z}{|f|})^{\frac{1}{2}}$ is the thickness of the Ekman layer and 
	the coefficients $k_1,\dots,k_4  \in \R$ are given in \eqref{eq:coefficients_2}.
\end{rem}

%

We conclude the following result which emphasizes the special role of the Ekman spiral in the regime $\alpha < 0$.

\begin{cor}[Uniqueness of equilibria]\label{cor:unique_2}
	Assume $\alpha < 0$. Then the Ekman spiral $(v_{\rE},0,\pi_{\rE})$ with pressure $\pi_{\rE}(x,y,z) = - \rho_0 g \cdot z$ is the unique steady-state with Coriolis force \eqref{eq:pe_2} and wind-driven boundary conditions \eqref{eq:wind bc_2} and geostrophic boundary conditions \eqref{eq:bottom_2}, i.e. the unique solution $(u,\pi)=(v,w,\pi)$ with $v \in \rH^1(\Omega)$ to the equations
	\begin{equation}
		\left\{
		\begin{aligned}
			- \nu_{\rH} \DeltaH v - \nu_z \partial_z^2 v + u \cdot \nabla v + \frac{1}{\rho_0}\nablaH \pi + f v^\perp &= 0, &&\text{ in } \Omega, \\ 
			\partial_z \pi+ \rho_0 g &= 0, &&\text{ in } \Omega, \\ 
			\div u &= 0, &&\text{ in } \Omega, 
		\end{aligned}
		\right. 
		\label{eq:equilibrium_2}
	\end{equation}
	and the boundary conditions \eqref{eq:wind bc_2}--\eqref{eq:w bc_2}.
\end{cor}
It remains an open question whether the Ekman spiral is the unique equilibrium when $\alpha \geq 0$. Moreover, it is unclear whether the smallness condition in \autoref{thm:main_2} in the transitional regime $\alpha \geq 0$ and $s(A) < 0$ can be removed.

\section{Ekman spiral}
\label{sec:Ekman_2}

In this section we study the Ekman spiral and establish some results which we need in the sequel. 
For convenience of the reader we show that the Ekman spiral satisfies the 3D incompressible primitive equations and Coriolis force \eqref{eq:pe_2} with \eqref{eq:wind bc_2} and \eqref{eq:bottom_2}.

\begin{prop}[Ekman spiral is an equilibrium]\label{prop:Ekman_2}
	The Ekman spiral \\ $(u_{\rE},\pi_{\rE}) := (v_{\rE},0,\pi_{\rE})$ given by \eqref{eq:ekman_2} and $\pi_{\rE}(x,y,z) := - \rho_0 g \cdot z$ satisfies the primitive equations \eqref{eq:pe_2} with wind-driven boundary conditions \eqref{eq:wind bc_2} and geostrophic boundary conditions \eqref{eq:bottom_2}. 
\end{prop}
\begin{proof}
	First, we remark by the explicit formulae from \autoref{sec:ekman explicit_2} that $v_{\rE} \in \rH^1(\Omega)$.  
	Note that the Ekman spiral is time-independent, i.e. $\partial_t v_{\rE} = 0$, and independent of $x,y$. The latter one
	implies $\nablaH v_{\rE} = 0$ and together with $w_{\rE} =0$ that 
	\begin{equation*}
		u_{\rE} \cdot \nabla v_{\rE}
		= v_{\rE} \cdot \nablaH v_{\rE} + w_{\rE} \partial_z v_{\rE}
		= 0 
	\end{equation*}
	and
	\begin{equation*}
		\div u_{\rE} = \divH v_{\rE} + \partial_z w_{\rE}
		= 0 .
	\end{equation*}
	Finally, we have $\nablaH \pi_{\rE} = 0$. 
	Using these identities and the fact that $v_{\rE}$ solves \eqref{eq:ekman_2} we obtain
	\begin{equation*}
		\partial_t v_{\rE} - \nu_{\rH} \DeltaH v_{\rE} - \nu_z \partial_z^2 v_{\rE} + u_{\rE} \cdot \nabla v_{\rE} + \frac{1}{\rho_0}\nablaH \pi_{\rE} + f v_{\rE}^\perp
		= - \nu_z \partial_z^2 v_{\rE} + f v_{\rE}^\perp
		= 0 . \qedhere 
	\end{equation*}
\end{proof}

Next, we prove \autoref{cor:unique_2}
under the assumption that \autoref{thm:main_2} holds.

\begin{proof}[Proof of \autoref{cor:unique_2}.]
	By \autoref{prop:Ekman_2} we know that $v_{\rE} \in \rH^1(\Omega)$ is an equilibrium. 	
	Now, let $v_* \in \rH^1(\Omega)$ be an equilibrium of \eqref{eq:pe_2}. It is in particular a solution of \eqref{eq:pe_2} with initial data $v(0) = v_\ast$. From \autoref{thm:main_2} and the fact that $v^\ast$ is time-independent we conclude
	\begin{equation*}
		v_\ast = \lim_{t\to+\infty}v^\ast(t) = v_{\rE} .
		\qedhere 
	\end{equation*}
\end{proof}

%

\section{Local wellposedness and Stability of the Ekman Spiral}
\label{sec:local_2}

In this section we introduce the functional analytic set-up to treat the primitive equations \eqref{eq:pe_2} with \eqref{eq:wind bc_2} and \eqref{eq:bottom_2}. Moreover, we provide necessary information about the linearization and prove local well-posedness of the system and finally asymptotic stability of the Ekman spiral. 

\smallskip 

We consider the differences $v_d := v - v_{\rE}$, $w_d := w$ and $\pi_d := \pi - \pi_{\rE}$. Using
\begin{align*}
	u_d \cdot \nabla v_d 
	+ v_{\rE} \cdot \nablaH v_d 
	+ w_d \partial_z v_{\rE} 
	&= (u-u_E) \cdot \nabla (v-v_{\rE}) 
	\\
	&+ v_{\rE} \cdot \nablaH (v-v_{\rE}) 
	+ w \partial_z v_{\rE} \\
	&= (v-v_{\rE}) \cdot \nablaH (v-v_{\rE})
	+ w \partial_z (v-v_{\rE}) \\
	&+ v_{\rE} \cdot \nablaH (v-v_{\rE}) 
	+ w \partial_z v_{\rE} \\
	&= v \cdot \nablaH (v-v_{\rE})
	+ w \partial_z v 
	= u \cdot \nabla v .
\end{align*}
we obtain from
\eqref{eq:pe_2}
the transformed system
\begin{equation}
	\left\{
	\begin{aligned}
		\partial_t v_d - \nu_{\rH} \DeltaH v_d - \nu_z \partial_z^2 v_d + u_d \cdot \nabla v_d \phantom{aaa} & && \\
		+ v_{\rE} \cdot \nablaH v_d 
		+ w_d \dz v_{\rE} + \frac{1}{\rho_0}\nablaH \pi_d + f v_d^\perp &= 0,
		&&\text{ in } (0,T) \times \Omega, \\
		\partial_z \pi_d &= 0, &&\text{ in } (0,T) \times \Omega, \\
		\div u_d &= 0, &&\text{ in } (0,T) \times \Omega, \\
		w_d &= 0, &&\text{ on } (0,T) \times \Gamma, \\
		\partial_z v_d &= 0 &&\text{ on } (0,T) \times \Gamma_u,\\
		v_d &= 0 &&\text{ on } (0,T) \times \Gamma_b, \\
		v_d(0) &= v_0-v_{\rE}, &&\text{ in } \Omega.
	\end{aligned}
	\right. 
	\label{eq:difference_2}
\end{equation}
We call these equations the \emph{Ekman-primitive equations}.
Next we introduce the spaces
\begin{equation*}
	X^k_0 := \rH^k(\Omega) \cap \rL^2_{\sigmabar}(\Omega)
\end{equation*}
and the
\emph{hydrostatic Ekman-Stokes operator} $A_k \colon X_1^k \subset X_0^k \to X_0^k$ by
\begin{equation}
	\begin{aligned} 
		A_k v_d &:= \P ( \nu_{\rH} \DeltaH v +\nu_z \partial_z^2 v_d)
		- \P (v_{\rE} \cdot \nablaH v_d 
		+ w_d \dz v_{\rE}) 
		- \P f v_d^\perp,
		\\
		X_1^k &:= \{ v \in \rH^{k+2}(\Omega) \cap \rL^2_{\sigmabar}(\Omega) \colon \partial_z v_d|_{\Gamma_u} = 0, \quad v_d|_{\Gamma_b} = 0 \} .
	\end{aligned} 
	\label{eq:A_2}
\end{equation} 
Note that $A_0$ coincides with the hydrostatic Ekman-Stokes operator $A$ from \autoref{sec:main_2}. Therefore this extends the definition of the hydrostatic Ekman-Stokes operator to Sobolev spaces. 
Note further that $X^k_0 \not = X_{k} := D(A^k)$\footnote{For an operator $A$ with domain $D(A)$ on a Banach space $X$ we define iteratively $D(A^k) := \{ x \in D(A^{k-1}) \colon A x \in D(A^{k-1}) \}$, see \cite[Section II.5]{EN:00_2}.} due to the presence of boundary conditions. 
Further, we define the bilinearity
\begin{equation*}
	F(v_d,v_d') := \P (v_d \cdot \nablaH v_d' + w(v_d) \dz v_d' ) 
\end{equation*}
with
\begin{equation}
	w(v_d)(x,y,z)
	= - \int_{-h}^z \divH v_d(x,y,\xi) \, \d \xi .
\end{equation}
Now \eqref{eq:difference_2} can be reformulated as the abstract Cauchy problem
\begin{equation}
	\dt v_d - A v_d = F(v,v), \qquad v_d(0) = v_0-v_{\rE} 
	\label{eq:acp_2}
\end{equation}
on $X_0$. 
Further, we introduce the maximal regularity spaces by
\begin{equation*}
	\E_1^k(T') := \rH^{1}(0,T';X_0^k) \cap \rL^2(0,T';X_1^k), \enspace \text{ and } \enspace 
	\E_1^k := \rH^{1}(\R_+;X_0^k) \cap \rL^2(\R_+;X_1^k). 
\end{equation*}
For simplicity of the notation we suppress the superscript $^k$ if $k =0$. 

\smallskip 

We denote the complex interpolation spaces  by $X_\theta := [X_0,X_1]_\theta$ for $\theta \in (0,1)$, where $[\cdot,\cdot]_\theta$ is the complex interpolation functor, and by $X_{\theta,p} := (X_0,X_1)_{\theta,p}$ for $\theta \in (0,1)$ and $p \in (1,\infty)$, where $(\cdot,\cdot)_{\theta,p}$ is the real interpolation functor
We recall from \cite{HHK:16_2} and \cite{GGHHK:20_2} the real and complex interpolation spaces for $p = q = 2$. 

\begin{lem}[Interpolation spaces]\label{lem:interpolation spaces_2}
	For $\theta \in (0,1)$ the complex interpolation spaces $X_\theta$ coincide with the real interpolation spaces $X_{\theta,2} := (X_0,X_1)_{\theta,2}$. They are explicitly given by
	\begin{equation*}
		X_\theta = X_{\theta,2}
		=  
		\begin{cases}
			\left\{ v_d \in \rH^{2\theta}(\Omega) \cap \rL^2_{\sigmabar}(\Omega) \colon \partial_z v_d|_{\Gamma_u} = 0 \quad  v_d|_{\Gamma_b} = 0 \right\}, \hspace{-1 em}&\text{ if } \theta \in (\nicefrac{3}{4},1), \\
			\left\{ v_d \in \rH^{2\theta}(\Omega) \cap \rL^2_{\sigmabar}(\Omega) \colon v_d|_{\Gamma_b} = 0 \right\}, &\text{ if } \theta \in (\nicefrac{1}{4},\nicefrac{3}{4}), \\
			\qquad \qquad \rH^{2\theta}(\Omega) \cap \rL^2_{\sigmabar}(\Omega), 
			&\text{ if } \theta \in (0,\nicefrac{1}{4}).
		\end{cases}
	\end{equation*} 
\end{lem}

In the next proposition we establish useful properties of the hydrostatic Ekman-Stokes operator. 

\begin{prop}[Bounded $\mathcal{H}^\infty$-calculus on $\rH^k$]\label{prop:mr_2}
	Consider $k \in \N_0$. 
	Then the operator $A_k+\omega$ for $\omega > s(A_k)$ admits a bounded $\mathcal{H}^\infty$-calculus on $X^k_0 = \rH^k(\Omega) \cap \rL^2_{\sigmabar}(\Omega)$
	and has compact resolvent on $X^k_0$. Hence, $A_k+\omega$ has maximal regularity, generates an analytic and compact semigroup on $X^k_0$ and its spectrum is independent of $k$ and consists only on discrete eigenvalues. If $k > 0$ the semigroup is not strongly continuous.  
\end{prop}
\begin{proof}
	We first show that the operator $A_k+\omega$ admits bounded $\mathcal{H}^\infty$-calculus on $X_0^k$ for $\omega > s(A_k)$. 
	In order to prove this we split the operator
	\begin{equation*}
		A_k = B_k + Q,
	\end{equation*}
	into the main part 
	\begin{equation*}
		B_k v_d := \P ( \nu_{\rH} \DeltaH v_d +\nu_z \partial_z^2 v_d),
		\quad 
		D(A_0) := X_1^k ,
	\end{equation*}
	and a perturbation coming from the linearization around the Ekman layer 
	\begin{equation*}
		Q v_d := 	- \P (v_{\rE} \cdot \nablaH v_d 
		+ w_d \dz v_{\rE}) 
		- \P f v_d^\perp .
	\end{equation*}
	Using the boundedness of the hydrostatic Helmholtz projection and the fact that $v_{\rE}$ is smooth, we see that
	\begin{equation*}
		Q \colon \rH^{k+1}(\Omega) \to \rH^k(\Omega)
	\end{equation*} 
	is bounded. Hence, $Q$ is relatively $(-B_k)^{\frac{1}{2}}$-bounded. Using \cite[Corollary 3.3.15]{PS:16_2} the problem reduces to the bounded $\mathcal{H}^\infty$-calculus of $B_k$ on $X_0^k$. In order to show this we follow the proof of \cite{GGHHK:17_2}.
	As shown in \cite{GGHHK:17_2}, by writing the pressure in terms of $v$, the resolvent problem
	\begin{equation*}
		\lambda v_d - \tilde{\Delta}v_d + \nablaH \pi_d = \lambda v_d - A_0 v_d = f 
	\end{equation*}
	for $\tilde{\Delta} \colon D(\tilde{\Delta}) \subset \rH^k(\Omega) \to \rH^k(\Omega)$ given by
	\begin{equation*}
		\tilde{\Delta} v_d := \nu_{\H}\Delta_{\H} v_d+ \nu_z \dz^2 v_d,
		\qquad D(\tilde{\Delta}) := \bigl\{ v_d \in \rH^{k+2}(\Omega) \colon \partial_z v_d|_{\Gamma_u} = 0, \ v_d|_{\Gamma_b} = 0  \bigr\} 
	\end{equation*}
	is equivalent to 
	\begin{equation*}
		\lambda v_d - \tilde{\Delta}v + C v_d = \sP f 
	\end{equation*}
	with 
	\begin{equation*}
		C v_d := -\nablaH \DeltaH^{-1} \divH \dz v_d\big|_{\Gamma_b} . 
	\end{equation*}
	Note that $C$ is a pseudo-differential operator of order $\nicefrac{3}{2}$ and hence
	relatively $(-\tilde{\Delta})^{\frac{3}{4}}$-bounded.
	Taking the vertical average and applying the horizontal divergence shows $\divH \vbar = 0$. Using the fact that $B$ is of lower order $\nicefrac{3}{2}$, elliptic regularity theory implies for $f \in \rH^k(\Omega)$ that $v \in \rH^k(\Omega)$. Hence, $R(\lambda,\tilde{\Delta}+C)$ maps $\rH^k(\Omega) \cap \rL^2_{\sigmabar}(\Omega)$ to itself and
	\begin{equation*}
		R(\lambda,B_k) = R(\lambda,\tilde{\Delta} +  C)\big|_{\rH^k(\Omega)\cap\rL^2_{\sigmabar}(\Omega)} .
	\end{equation*}
	Using \cite[Corollary 3.3.15]{PS:16_2} the problem reduces to the bounded $\mathcal{H}^\infty$-calculus of $\tilde{\Delta}$ on $\rH^k(\Omega)$.
	By rescaling in $x,y$ or $z$ it suffices to consider the case of the classic Laplacian $\Delta$ instead of $\tilde{\Delta}$.
	As in \cite[Lemma 4.1]{GGHHK:17_2} using extension and restriction arguments this can be reduced to the problem of the bounded $\mathcal{H}^\infty$-calculus of $\Delta_{\T^3,k}$ on $\rH^k(\T^3)$.
	For $k = 0$ this is shown in \cite[Lemma 4.1]{GGHHK:17_2} by using the parameter ellipticity of the symbol. For $k \in \N$ we note that
	\begin{equation*}
		D(\Delta_{\T^3,k}) = \rH^{k+2}(\T^3) = D(\Delta_{\T^3}^k) 
	\end{equation*}
	is an interpolation space of $\Delta_{\T^3}$
	and $\Delta_{\T^3,k} = \Delta_{\T^3}\big|_{\rH^k(\T^3)}$. By \cite[Corollary II.5.3]{EN:00_2} the operators $\Delta_{\T^3,k}$ are similar for all $k \in \N_0$. Hence, the bounded $\mathcal{H}^\infty$-calculus of $\Delta_{\T^3,k}$ follows from the bounded $\mathcal{H}^\infty$-calculus of $\Delta_{\T^3,0}=\Delta_{\T^3,0}$. 
	
	\medskip 
	
	Since $X_1^k \stackrel{c}{\hookrightarrow} X_0^k$ it follows that $A_k$ has compact resolvent and hence its spectrum consists only on eigenvalues, i.e. $\sigma(A_k) = \sigma_p(A_k)$. 
	Since eigenfunctions are smooth, the spectrum is $k$-independent. 
\end{proof} 

We need the following improvement of the bounds from \cite[Corollary 6.2]{GGHHK:20_2}.

\begin{lem}[Bilinear bounds]\label{lem:bilinear bounds_2}
	Let $k \in \N$. Then 
	\begin{equation*}
		\| F (v_d,v_d) \|_{\rH^{k}(\Omega)}
		\leq C \cdot 
		\| v_d \|_{\rH^{k+2}(\Omega)}^{\frac{1}{2}} \cdot \| v_d \|_{\rH^{k+1}(\Omega)} \cdot 
		\| v_d \|_{\rH^{k}(\Omega)}^{\frac{1}{2}} .
	\end{equation*}
	Moreover, for $k = 0$ we obtain
	\begin{equation*}
		\| F (v_d,v_d) \|_{\rL^2(\Omega)}  
		\leq C \cdot 
		\| v_d \|_{\rH^{\frac{3}{2}}(\Omega)}^2 .
	\end{equation*}
\end{lem}
\begin{proof}
	The result for $k = 0$ is proven in \cite[Lemma 5.1]{HK:16_2}.
	For $k = 1$ we note that
	\begin{equation*}
		\nabla F(v_d,v_d)
		= \nabla v_d \cdot \nablaH v_d + v_d \cdot \nablaH \nabla v_d + \nabla w_d \partial_z v_d  + w_d \partial_z \nabla v_d . 
	\end{equation*}
	We estimate the terms separately. 
	Using $\rH^2(\Omega) \hookrightarrow \rH^{1,4}(\Omega)$ and the interpolation inequality, we obtain for the first term
	\begin{align*}
		\| \nabla v_d \cdot \nablaH v\|_{\rL^2(\Omega)}
		&\leq C \cdot \| v_d \|_{\rH^{1,4}(\Omega)}^2
		\leq C \cdot \| v_d \|_{\rH^{2}(\Omega)}^2 \\
		&\leq C \cdot \| v_d \|_{\rH^{3}(\Omega)}^{\frac{1}{2}} \cdot 
		\| v_d \|_{\rH^{2}(\Omega)} \cdot 
		\| v_d \|_{\rH^{1}(\Omega)}^{\frac{1}{2}} .
	\end{align*}
	For the second term we obtain by Agmon's inequality
	\begin{equation*}
		\| v_d \cdot \nablaH \nabla v_d \|_{\rL^2(\Omega)}
		\leq C \cdot \| v_d \|_{\rL^\infty(\Omega)} \cdot \| v_d \|_{\rH^2(\Omega)}
		\leq C \cdot \| v_d \|_{\rH^2(\Omega)}^{\frac{3}{2}} \cdot \| v_d \|_{\rH^1(\Omega)}^{\frac{1}{2}} .
	\end{equation*}
	For the third term we obtain from the embedding $\rH^{\frac{1}{2}}(\T^2) \hookrightarrow \rL^4(\T^2)$ and the interpolation inequality
	\begin{align*}
		\| \nabla w_d \partial_z v_d \|_{\rL^2(\Omega)}
		&\leq \| \nabla w_d \|_{\rL^\infty_z \rL^4_{xy}(\Omega)} \cdot \| \partial_z v_d \|_{\rL^2_z \rL^4_{xy}(\Omega)} \\
		&\leq C \cdot \| \nabla \divH v_d \|_{\rL^2_z \rL^4_{xy}(\Omega)} \cdot \| \partial_z v_d \|_{\rL^2_z \rL^4_{xy}(\Omega)} \\
		&\leq C \cdot \| \nabla \divH v_d \|_{\rL^2_z \rH^{\frac{1}{2}}_{xy}(\Omega)} \cdot \| \partial_z v_d \|_{\rL^2_z \rH^{\frac{1}{2}}_{xy}(\Omega)} \\
		&\leq C \cdot \| \nabla \divH v_d \|_{\rH^{\frac{1}{2}}(\Omega)} \cdot \| \partial_z v_d \|_{\rH^{\frac{1}{2}}(\Omega)} \\
		&\leq C \cdot \| v_d \|_{\rH^{\frac{5}{2}}(\Omega)} \cdot 
		\| v_d \|_{\rH^{\frac{3}{2}}(\Omega)} \\
		&\leq  C \cdot 
		\| v_d \|_{\rH^{3}(\Omega)}^{\frac{1}{2}} \cdot \| v_d \|_{\rH^{2}(\Omega)} \cdot 
		\| v_d \|_{\rH^{1}(\Omega)}^{\frac{1}{2}} 
	\end{align*}
	Similarly we obtain for the fourth term 
	\begin{align*}
		\| w_d \nabla \partial_z v_d \|_{\rL^2(\Omega)}
		&\leq \| w_d \|_{\rL^\infty_z \rL^4_{xy}(\Omega)} \cdot \| \nabla \partial_z v_d \|_{\rL^2_z \rL^4_{xy}(\Omega)} \\
		&\leq C \cdot \| \divH v_d \|_{\rL^2_z \rL^4_{xy}(\Omega)} \cdot \| \nabla \partial_z v_d \|_{\rL^2_z \rL^4_{xy}(\Omega)} \\
		&\leq C \cdot \| \divH v_d \|_{\rL^2_z \rH^{\frac{1}{2}}_{xy}(\Omega)} \cdot \| \nabla \partial_z v_d \|_{\rL^2_z \rH^{\frac{1}{2}}_{xy}(\Omega)} \\
		&\leq C \cdot \| \divH v_d \|_{\rH^{\frac{1}{2}}(\Omega)} \cdot \| \nabla \partial_z v_d \|_{\rH^{\frac{1}{2}}(\Omega)} \\
		&\leq C \cdot \| v_d \|_{\rH^{\frac{3}{2}}(\Omega)} \cdot \| v_d \|_{\rH^{\frac{5}{2}}(\Omega)} \\
		&\leq  C \cdot \| v_d \|_{\rH^{3}(\Omega)}^{\frac{1}{2}} \cdot \| v_d \|_{\rH^{2}(\Omega)} \cdot \| v_d \|_{\rH^{1}(\Omega)}^{\frac{1}{2}} .
	\end{align*}
	Finally, for $k \geq 2$ we have that $\rH^k(\Omega)$ is a Banach algebra. Hence using the interpolation inequality
	\begin{align*}
		\| F (v_d,v_d) \|_{\rH^{k}(\Omega)}
		&\leq C\cdot 
		\| v_d \|_{\rH^{k}(\Omega)} \cdot \| \nablaH v_d \|_{\rH^k(\Omega)}
		+ C \cdot \| w_d \|_{\rH^k(\Omega)} \cdot \| \partial_z v_d \|_{\rH^{k}(\Omega)} 
		\\
		&\leq C \cdot \| v_d \|_{\rH^{k+1}(\Omega)}^2
		\leq C \cdot \| v_d \|_{\rH^{k+2}(\Omega)}^{\frac{1}{2}} \cdot \| v_d \|_{\rH^{k+1}(\Omega)} \cdot \| v_d \|_{\rH^{k}(\Omega)}^{\frac{1}{2}} . 
	\end{align*}
\end{proof}

\begin{prop}[Local wellposedness]\label{prop:local_2}
	For every $v_0 \in \rH^{1}(\Omega) \cap \rL^2_{\sigmabar}(\Omega)$ satisfying \eqref{eq:bottom_2} there exists $T' > 0$ and a unique solution $v_{\textcolor{red}{d}} \in \E_1(T')$ to \eqref{eq:difference_2}. 
	This implies the existence of unique solution $v \in \rL^2(0,T';\rH^2(\Omega)) \cap \rH^1(0,T;\rL^2_{\sigmabar}(\Omega)) \hookrightarrow \mathrm{BUC}([0,T');\rH^{1}(\Omega)\cap \rL^2_{\sigmabar}(\Omega))$ to \eqref{eq:pe_2} with the boundary conditions \eqref{eq:wind bc_2}--\eqref{eq:w bc_2}.
\end{prop}
\begin{proof}
	Our proof follows the strategy from \cite{GGHHK:20_2} which verifies the assumption from \cite{PW:17_2}. Note that (H1) is trivial here, since we consider $F_1 = 0$. Furthermore, conditions (H2) and (H3) follow from \autoref{lem:bilinear bounds_2} with $\beta = \frac{3}{4}$.
	It remains to verify the maximal regularity and condition (S): 
	By \cite[Remark 1.1]{PW:17_2} it is sufficient to show that $A+\omega$ admits bounded $\mathcal{H}^\infty$-calculus on $X_0$ which is proven in \autoref{prop:mr_2}. 
	Now the claim follows from \cite[Theorem 1.2]{PW:17_2}. 
	Using the decomposition $v = v_d + v_{\rE}$ and	$v_d \in \E_1(T')$ and $v_{\rE} \in \rC^\infty(-h,0)$ we conclude the desired regularity of $v$.
	For the boundary conditions, we obtain  
	\begin{equation*}
		\partial_z v = \partial_z v_{\rE} + \partial_z v_d = \partial_z v_{\rE} = \tau \text{ on } \Gau, \qquad \text{ and } \qquad
		v = v_{\rE} + v_d = v_{\rE} = v_g \text{ on } \Gab . \qedhere
	\end{equation*} 
\end{proof}

Finally, we conclude the stability of the Ekman spiral in the primitive equations if $s(A) < 0$. Although the following result is a direct consequence of the principle of linearized stability and could be concluded from \cite[Theorem 5.3.1 \& Remark 5.3.2]{PS:16_2} we provide a direct proof here for convince. 

\begin{prop}[Stability of the Ekman spiral]\label{cor:small data_2}
	If $s(A) < 0$ holds
	then there exists $r > 0$ such that for every $v_0 \in \rH^{1}(\Omega) \cap \rL^2_{\sigmabar}(\Omega)$ satisfying \eqref{eq:bottom_2}
	and $\| v_0-v_{\rE} \|_{\rH^1(\Omega)} < r$ there exists 
	a unique global solution $v_d \in \E_1$ to \eqref{eq:difference_2} and $\| v_d(t) \|_{\rH^1(\Omega)} \leq C e^{\omega t}$ for all $\omega < s(A) < 0$, where $s(A)$ is the spectral bound of $A$.  
\end{prop}
\begin{proof}
	Note that we have seen in the proof of \autoref{prop:local_2} that the assumptions from \cite{PWS:18_2} are satisfied.
	Since $s(A) < 0$, we obtain from \autoref{prop:mr_2} that the semigroup generated by $A$ is exponentially stable.  
	Recall that $A$ admits maximal $\rL^p$-regularity for all $p \in (1,+\infty)$ and $\beta = \frac{3}{4}$.
	For $p \in (4,\infty)$ we have $\frac{1}{p} < 1-\beta$. 
	It follows from \cite[Corollary 2.2]{PWS:18_2} that the solution $v$ to \eqref{eq:difference_2} decays exponentially in $X_{1-\nicefrac{1}{p},p}$,
	if $\| v_0 - v_{\rE} \|_{X_{\mu_c-\frac{1}{p},p}} < r_0$ with $\mu_c := 2 \beta - 1+\frac{1}{p} = \frac{1}{2}+ \frac{1}{p}$.
	By \cite[Lemma 2.1]{GGHHK:20_2} we know that
	$X_{1-\nicefrac{1}{p},p} = \{ v_d \in \rB^{2-\nicefrac{2}{p}}_{2p}(\Omega) \cap \rL^2_{\sigmabar}(\Omega) \colon \partial_z v_d|_{\Gau} = 0, \quad v_d|_{\Gab} = 0 \}$
	and 
	$X_{\mu_c-\nicefrac{1}{p},p} =
	X_{\frac{1}{2},p}
	= \{ v_d \in \rB^{1}_{2p}(\Omega) \cap \rL^2_{\sigmabar}(\Omega)\colon v_d|_{\Gab} = 0 \}$, where $\rB^{s}_{2pp}(\Omega)$ denotes the Besov space of order $s$.
	From the embedding 
	\begin{equation*}
		\rH^1(\Omega)
		= \rB^1_{22}(\Omega)
		\hookrightarrow \rB^{1}_{2p}(\Omega)
	\end{equation*}
	we obtain
	\begin{equation*}
		\| v_0 - v_{\rE} \|_{\rB^{1-\nicefrac{2}{p}}_{2pp}(\Omega)}
		\leq 
		C \cdot 
		\| v_0 - v_{\rE} \|_{\rH^1(\Omega)}
		\leq C \cdot r .
	\end{equation*}
	Thus choosing $r>0$ sufficiently small such that $r_0 := C \cdot r$ is sufficiently small we obtain exponential decay in 	$X_{1-\nicefrac{1}{p},p} = \rB^{2-\nicefrac{2}{p}}_{2p}(\Omega) \cap \rL^2_{\sigmabar}(\Omega)$.
	Using the embedding 
	\begin{equation*}
		\rB^{2-\nicefrac{2}{p}}_{2p}(\Omega) \hookrightarrow \rH^{1}(\Omega)
	\end{equation*}
	we obtain 
	\begin{equation*} 
		\| v_d(t) \|_{\rH^1(\Omega)}
		\leq \| v_d(t) \|_{\rB^{2-\nicefrac{2}{p}}_{2p}(\Omega)} \leq C e^{\omega t}
	\end{equation*}
	for $\omega < \omega_0 = s(A) < 0$, where we have used \cite[Corollary IV.3.12]{EN:00_2}.
\end{proof}

We finish this section with the instability of the Ekman spiral in the primitive equations if $s(A) > 0$. 

\begin{prop}[Instability of the Ekman spiral]\label{prop:instability_2}
	If $s(A) > 0$ holds then there exists $\varepsilon >0$ such that for any $\delta > 0$ there exists $v_0 \in \rH^{1}(\Omega) \cap \rL^2_{\sigmabar}(\Omega)$ satisfying \eqref{eq:bottom_2} with $\| v_0 - v_{\rE} \|_{\rH^1(\Omega)} < \delta$ such that there is a finite time $t_0 > 0$ with
	\begin{equation*}
		\| v_d(t_0) \|_{\rH^1(\Omega)} \geq \varepsilon . 
	\end{equation*}
\end{prop}
\begin{proof}
	We verify the conditions from \cite[Theorem 5.4.1]{PS:16_2}. The maximal regularity of $A$ has been shown in \autoref{prop:mr_2}. Moreover, we have shown that the spectrum is discrete and consists purely of eigenvalues. By assumption we know that the largest eigenvalue $\Re \lambda_1 = s(A) > 0$. Since, the spectrum is discrete, we know that there exists $\eta > 0$ such that $\Re\lambda_1 - \eta + i \R \in \rho(A)$. 
	Now, the claim follows from  \cite[Theorem~5.4.1]{PS:16_2}.
\end{proof} 

\section{Global wellposedness and Convergence}
\label{sec:global_2}

In this section we show the global well-posedness of the primitive equations \eqref{eq:pe_2} with \eqref{eq:wind bc_2} and \eqref{eq:bottom_2} and conclude together with our stability result \autoref{cor:small data_2} the convergence of the solutions to the Ekman spiral as $t \to + \infty$. 
In order to show the global well-posedness so we first describe the blow-up scenario which we will rule out in the sequel.
From \cite[Theorem 2.4]{PWS:18_2} we obtain the following blow-up criterion. 

\begin{prop}[Blow-up criterion]
	\label{cor:blow-up_2}
	Let $v_d$ be the local solution from \autoref{prop:local_2} with maximal interval of existence $[0,T_+)$. If $T_+ < +\infty$ then 
	\begin{equation*}
		\limsup_{T \to T_+} \int_0^T \|v_d(s)\|_{\rH^{2}(\Omega)}^2 \d s = +\infty . 
	\end{equation*}
\end{prop}
\begin{proof}
	Note that we have seen in the proof of \autoref{prop:local_2} that the assumptions from \cite{PWS:18_2} are satisfied.
	Moreover, we conclude from the proof of \autoref{prop:local_2} that the critical weight is given by
	$\mu_c := 2 \beta - 1+\frac{1}{2} = \frac{1}{2}+ \frac{1}{2}=1$. This means that the critical (complex) interpolation space is $X_{\mu_c} = X_1$. Now the result follows by \cite[Theorem 2.4]{PWS:18_2}. 
\end{proof}

Now, we show a-priori estimates. The general strategy follows the strategy from \cite{HK:16_2, GGHHK:20_2}. 
We start with energy estimates. 
We first need a G\aa rding inequality for the hydrostatic Ekman-Stokes operator. 

\begin{lem}[G\aa rding inequality]\label{lem:garding inequality_2}
	The hydrostatic Ekman-Stokes operator $A$ on $X_0$ satisfies 
	\begin{equation*}
		- \Re \langle A v_d, v_d \rangle 
		\geq c_0 \cdot \| \nabla v_d \|_{\rL^2(\Omega)}^2 - K \cdot \| v_d \|_{\rL^2(\Omega)}^2
	\end{equation*}
	for all $v_d \in D(A)$
	with constants $c_0 > 0$ and $K > 0$ depending on $\| \partial_z v_{\rE} \|_{\rL^\infty(-h,0)}$. 
\end{lem}
\begin{proof}
	We denote the complex conjugate of $v$ by $v^\ast$. 
	Using 
	\begin{equation*}
		\int_{\Omega} \nablaH \pi v_d^\ast 
		= \int_{\T^2} \nablaH \pi \int_0^h v_d^\ast 
		= h \int_{\T^2} \nabla \pi \vbar_d^\ast 
		= h \int_{\T^2} \pi_d \cdot  (\divH\vbar_d)^\ast = 0,
	\end{equation*}
	and 
	\begin{equation*}
		v_d^\perp v_d^\ast + (v_d^\ast)^\perp v_d 
		= -v_2 v_1^\ast + v_1 v_2^\ast-v_2^\ast v_1 + v_1^\ast v_2 = 0 
	\end{equation*}
	we obtain by 
	multiplying by $v_d^\ast$  
	and integration by parts
	\begin{align*}
		2\Re \langle A v_d, v_d \rangle 
		&= 2\int_{\Omega} B v_d v_d^* + \int_{\Omega} Q v_d v_d^* \\
		&= -2\nu_{\rH} \| \nablaH v_d \|_{\rL^2(\Omega)}^2
		- 2\nu_z \| \dz v_d \|_{\rL^2(\Omega)}^2
		+\Re \int_{\Omega} v_{\rE} \cdot \nablaH v_d v_d^\ast \\
		&+ \Re\int_{\Omega} v_{\rE} \cdot \nablaH v_d^\ast v_d 
		+ \Re \int_{\Omega} w_d \dz v_{\rE} v_d^\ast
		+ \Re \int_{\Omega} w_d^\ast \dz v_{\rE} v_d .
	\end{align*}
	Using integration by parts for
	the third term on the right-hand side yields
	\begin{equation*}
		\int_{\Omega} v_{\rE} \cdot \nablaH v_d v_d^\ast 
		= - \int_{\Omega} \divH v_{\rE} |v_d|^2 -\int_{\Omega} v_{\rE} \cdot \nablaH v_d^\ast v_d
		= -\int_{\Omega} v_{\rE} \cdot \nablaH v_d^\ast v_d
	\end{equation*}
	and hence it cancels with the forth term. By \autoref{rem:w_2} and Jensen's inequality we obtain
	\begin{equation}
		\label{eq:w jensen_2}
		\begin{aligned} 
			\| w_d \|_{\rL^2(\Omega)}^2
			&= \int_{-h}^0 \int_{\T^2} \left| \int_{-h}^z \divH v_d \right|^2
			\leq \int_{-h}^0 \int_{\T^2} \left( \int_{-h}^z |\divH v_d| \right)^2 \\
			&\leq \int_{-h}^0 \int_{\T^2} h^2 \cdot \left( \frac{1}{h} \int_{-h}^0 |\divH v_d| \right)^2 
			\leq h \int_{\T^2} h^2 \cdot \frac{1}{h} \int_{-h}^0 |\divH v_d|^2  \\
			&\leq 2 h^2 \int_{\T^2} \int_{-h}^0 |\nablaH v_d|^2 = 2 h^2 \cdot \| \nablaH v_d \|_{\rL^2(\Omega)}^2 .
		\end{aligned}
	\end{equation}
	Therefore, we obtain for the fifth and sixth term we obtain by Cauchy-Schwarz inequality
	\begin{align*}
		\int_{\Omega} |w_d \dz v_{\rE} v_d^\ast|
		&\leq \| \dz v_{\rE} \|_{\rL^\infty(-h,0)} \cdot \biggl( \frac{\varepsilon}{2} \cdot \| w_d \|_{\rL^2(\Omega)}^2
		+ \frac{1}{2\varepsilon} \cdot \| v_d \|_{\rL^2(\Omega)}^2 \biggr)\\
		&\leq \| \dz v_{\rE} \|_{\rL^\infty(-h,0)} \cdot \biggl( h^2 \varepsilon \cdot \| \nablaH v_d \|_{\rL^2(\Omega)}^2
		+ \frac{1}{2\varepsilon} \cdot 
		\| v_d \|_{\rL^2(\Omega)}^2 \biggr) .
	\end{align*}
	Choosing $\varepsilon = \tfrac{\nu_{\rH}}{16 h^2 \| \partial_z v_{\rE} \|_{\rL^\infty(-h,0)}}$ we obtain  
	\begin{align*}
		2\Re \langle A v_d, v_d \rangle 
		&= -2\nu_{\rH} \| \nablaH v_d \|_{\rL^2(\Omega)}^2
		- 2\nu_z \| \dz v_d \|_{\rL^2(\Omega)}^2 \\
		&+ \Re \int_{\Omega} w_d \dz v_{\rE} v_d^\ast
		+ \Re \int_{\Omega} w_d^\ast \dz v_{\rE} v_d \\
		&\leq -2\nu_{\rH} \| \nablaH v_d \|_{\rL^2(\Omega)}^2
		- 2\nu_z \| \dz v_d \|_{\rL^2(\Omega)}^2 \\
		&+ \nu_{\rH} \| \nablaH v_d \|_{\rL^2(\Omega)}^2
		+  K \cdot \| v_d \|_{\rL^2(\Omega)}^2 \\
		&= -\nu_{\rH} \| \nablaH v_d \|_{\rL^2(\Omega)}^2
		- 2\nu_z \| \dz v_d \|_{\rL^2(\Omega)}^2
		+ K \cdot \| v_d \|_{\rL^2(\Omega)}^2  
	\end{align*}
	and the desired inequality follows for $c_0 := \min\{ v_{\rH},2\nu_z \} > 0$.
\end{proof}

While Garding's inequality holds without any restrictions of the parameters, we obtain a stronger inequality assuming that $\alpha < 0$.

\begin{cor}[Coercivity]\label{cor:coercivity_2}
	Suppose $\alpha < 0$. Then there exists a constant $c_0>0$ such that
	\begin{equation*}
		- \Re \langle A v_d, v_d \rangle 
		\geq c_0 \cdot 
		\| \nabla v_d \|_{\rL^2(\Omega)}^2 
	\end{equation*}
	for all $v_d \in D(A)$.
\end{cor}
\begin{proof}
	By definition of $\alpha$ we obtain
	\begin{equation*}
		\Re \langle Av_d , v_d \rangle \leq \alpha \cdot \| v_d \|_{\rL^2(\Omega)}^2 . 
	\end{equation*}
	Using \autoref{lem:garding inequality_2} this implies
	\begin{align*}
		- \Re \langle A v_d, v_d \rangle 
		&= - \delta \cdot \Re \langle A v_d, v_d \rangle 
		- (1-\delta) \cdot \Re \langle A v_d, v_d \rangle \\
		&\geq \delta \cdot  c_0 \cdot \| \nabla v_d \|_{\rL^2(\Omega)}^2 - \delta \cdot K \cdot \| v_d \|_{\rL^2(\Omega)}^2 - (1-\delta) \cdot \alpha \cdot 
		\| v_d \|_{\rL^2(\Omega)}^2 
	\end{align*}
	for $\delta \in (0,1)$. Choosing $\delta \leq \frac{-\alpha}{K-\alpha}$ we obtain 
	\begin{equation*}
		- \delta \cdot K \| v_d \|_{\rL^2(\Omega)}^2 
		- (1-\delta) \cdot \alpha \cdot 
		\| v_d \|_{\rL^2(\Omega)}^2 
		= ( \delta \cdot (\alpha-K) - \alpha  ) \cdot \| v_d \|_{\rL^2(\Omega)}^2 
		\geq 0 
	\end{equation*}
	and the claim follows with $\tilde{c}_0 := c_0 \cdot \delta > 0$. 
\end{proof}

\begin{rem}
	Since for real-valued vector-fields $v_d,v_d' \in D(A)$ the expression $\langle Q v_d,v_d' \rangle$ is always real, the coercivity condition from \autoref{cor:coercivity_2}
	is equivalent to 
	\begin{equation*}
		- \langle A v_d, v_d \rangle 
		\geq c_0 \cdot 
		\| \nabla v_d \|_{\rL^2(\Omega)}^2 
	\end{equation*} 
	where only real-valued vector-fields are considered. The constant $c_0>0$ is identical. 
\end{rem}

\begin{lem}[Energy estimate]\label{lem:energy_2}
	Suppose $\alpha < 0$. 
	Let $v_d$ be a solution to \eqref{eq:difference_2}. Then
	\begin{equation*}
		\| v_d(t) \|_{\rL^2(\Omega)}^2
		+ \int_0^t \| \nabla v_d(s) \|_{\rL^2(\Omega)}^2 \d s \leq C \cdot \| v_0-v_{\rE} \|_{\rL^2(\Omega)}^2 .
	\end{equation*}
\end{lem}
\begin{proof} 	
	Testing \eqref{eq:difference_2}$_1$ with $v_d$ and using the cancellation law
	\begin{equation*}
		\langle \nablaH u_d \cdot \nabla v_d , v_d \rangle 
		=
		2 \int_{\Omega} u_d \cdot \nabla v_d \cdot v_d 
		= -\int_{\Omega} \div(u_d) \cdot |v_d|^2 = 0,
	\end{equation*}
	and
	\begin{equation*}
		\langle \nablaH \pi , v_d \rangle =
		\int_{\Omega} \nablaH \pi v_d 
		= \int_{\T^2} \nablaH \pi \int_0^h v_d
		= h \int_{\T^2} \nabla \pi \vbar_d
		= h \int_{\T^2} \pi \divH\vbar_d = 0,
	\end{equation*}
	implies 
	\begin{equation*}
		\frac{1}{2}\cdot \dt \| v_d \|_{\rL^2(\Omega)}^2
		- \langle A v_d , v_d \rangle \leq 0 .
	\end{equation*}
	Now, \autoref{cor:coercivity_2} implies 
	\begin{equation*}
		\frac{1}{2} \cdot \dt \| v_d \|_{\rL^2(\Omega)}^2
		+ c_0 \cdot \| \nabla v_d \|_{\rL^2(\Omega)}^2 \leq 0 .
	\end{equation*}
	and the desired estimate follows by integrating in time. 
\end{proof}

In the next step we establish an estimate for $\partial_z v$. 
\begin{lem}[Estimates for $\partial_z v_d$]\label{lem:estimate vz_2}
	Suppose $\alpha < 0$.  
	Then
	\begin{equation*}
		\begin{aligned}
			&\| \partial_z v_d(t) \|_{\rL^2(\Omega)}^2 
			+ \nu_{\rH} \int_0^t \| \nablaH \partial_z v_d \|_{\rL^2(\Omega)}^2
			+ \nu_z \int_0^t \| \partial_z^2 v_d \|_{\rL^2(\Omega)} \\ 
			\leq \ &\| \partial_z (v_0 - v_{\rE}) \|_{\rL^2(\Omega)}^2
			+ 
			C \cdot (1 + \| \partial_z v_{\rE} \|_{\rL^\infty(-h,0)} + \| \partial_z^2 v_{\rE} \|_{\rL^2(-h,0)}) \cdot \| v_0 - v_{\rE} \|_{\rL^2(\Omega)}^2 \\
			+ \ &\varepsilon_z \cdot \int_0^t \| \nablaH \pi \|_{\rL^2(\T^2)}^2
			+ C \cdot \int_0^t \| \nabla v_d\|_{\rL^2(\Omega)}^2 \\
			+ \ &C \cdot \int_0^t (\| \nablaH \vbar_d \|_{\rL^2(\Omega)} + \| \nablaH \vbar_d \|_{\rL^2(\Omega)}^2) \cdot \| \partial_z v \|_{\rL^2(\Omega)}^2 
			+ M_z \cdot \int_0^t \| |\vtilde_d| |\nabla\vtilde_d| \|_{\rL^2(\Omega)}^2  
		\end{aligned}
	\end{equation*}
	for $\varepsilon_z > 0$ and a fixed constant $M_z > 0$.
\end{lem}
\begin{proof}
	As in \cite[p. 2743]{KZ:07_2}, \cite[Step 2]{HK:16_2} or \cite[Step 3]{GGHHK:20_2} we 
	multiply \eqref{eq:difference_2} by $-\partial_z^2 v_d$ and integrate by parts
	\begin{equation*}
		\begin{aligned}
			&\frac{1}{2} \partial_t \| \partial_z v_d \|_{\rL^2(\Omega)}^2 
			+ \nu_{\rH} \| \nablaH \partial_z v_d \|_{\rL^2(\Omega)}^2
			+ \nu_z \| \partial_z^2 v_d \|_{\rL^2(\Omega)}
			\\
			= \ & 
			- \frac{1}{\rho_0} \int_{\T^2} (\nablaH \pi \cdot \partial_z v_d)\big|_{\Gamma_b} 
			- \int_\Omega (\partial_z v_d \cdot \nablaH v_d) \cdot \partial_z v_d 
			+ \int_\Omega \divH v_d \cdot \partial_z v_d \cdot \partial_z v_d \\
			\ &- \int_\Omega (\partial_z v_{\rE} \cdot \nablaH v_d) \cdot \partial_z v_d 
			+ \int_\Omega w_d \cdot \partial_z v_{\rE} \cdot \partial_z^2 v_d \\
			=: \ & \I + \II + \III + \IV + \V, 
		\end{aligned}
	\end{equation*}
	where have used 
	\begin{equation*}
		\int_\Omega f \partial_z v_d^\perp \cdot \partial_z^2 v_d
		= - \int_\Omega f \partial_z v_d^\perp \cdot \partial_z v_d
		= - \int_\Omega f (\partial_z v_d)^\perp \cdot  \partial_z v_d
		= 0 
	\end{equation*}
	and
	\begin{equation*}
		\int_\Omega (v_{\rE} \cdot \nablaH v_d) \cdot \partial_z^2 v_d 
		= 
		- \int_\Omega \partial_z v_{\rE} \cdot \nablaH v_d \cdot \partial_z v_d .
	\end{equation*}
	From \cite[Step 2, (6.7)]{HK:16_2} we know that
	\begin{align*}
		| \I |
		&\leq \frac{1}{\rho_0} \cdot\| \nablaH \pi \|_{\rL^2(\T^2)}
		\cdot \| \partial_z v_d\|_{\rL^{2}(\Omega)}^{\frac{1}{2}}
		\cdot \| \partial_z v_d\|_{\rH^{1}(\Omega)}^{\frac{1}{2}} \\
		&\leq \varepsilon_1 \| \nablaH \pi \|_{\rL^2(\T^2)}^2
		+ \frac{\varepsilon_2}{4 \rho_0 \cdot \varepsilon_1} \| \partial_z v_d\|_{\rH^{1}(\Omega)}^2
		+ C \cdot \| \partial_z v_d\|_{\rL^2(\Omega)}^2
	\end{align*}
	for $\varepsilon_1,\varepsilon_2 >0$,
	and
	\begin{align*}
		|\II| + |\III|
		&\leq
		C \cdot(\| \nablaH \vbar_d \|_{\rL^2(\Omega)} + \| \nablaH \vbar_d \|_{\rL^2(\Omega)}^2) \cdot \| \partial_z v \|_{\rL^2(\Omega)}^2 \\
		&+ M_z \cdot \| |\vtilde_d| |\nabla\vtilde_d| \|_{\rL^2(\Omega)}^2
		+ \varepsilon_3 \| \nabla \partial_z v_d \|_{\rL^2(\Omega)}^2 
	\end{align*}
	for $\varepsilon_3 > 0$.
	It remains to control the two terms involving the Ekman spiral. 
	For $\IV$ we obtain
	\begin{equation*}
		| \IV |
		\leq 
		\| \partial_z v_{\rE} \|_{\rL^\infty(-h,0)}
		\cdot \| \nabla v_d \|_{\rL^2(\Omega)}^2 .
	\end{equation*}
	Using \autoref{lem:energy_2} it follows
	\begin{equation*}
		\int_0^t |\IV|
		\leq \| \partial_z v_{\rE} \|_{\rL^\infty(-h,0)} 
		\cdot \int_0^t \| \nabla v_d \|_{\rL^2(\Omega)}^2
		\leq C \cdot \| \partial_z v_{\rE} \|_{\rL^\infty(-h,0)} \cdot \| v_0 - v_{\rE} \|_{\rL^2(\Omega)}^2 .
	\end{equation*}
	For \V \ we note that 
	\begin{align*}
		\int_\Omega w_d \cdot \partial_z v_{\rE} \cdot \partial_z^2 v_d
		&= - \int_\Omega \partial_z w_d \cdot \partial_z v_{\rE} \cdot \partial_z v_d
		- \int_\Omega  w_d \cdot \partial_z^2 v_{\rE} \cdot \partial_z v_d \\ 
		&= \int_\Omega \divH v_d \cdot \partial_z v_{\rE} \cdot \partial_z v_d
		- \int_\Omega  w_d \cdot \partial_z^2 v_{\rE} \cdot \partial_z v_d \\
		&=: \VI + \VII . 
	\end{align*}
	For the first term we obtain
	\begin{equation*}
		|\VI|
		\leq \| \partial_z v_{\rE} \|_{\rL^\infty(-h,0)}
		\cdot  \| \nabla v_d \|_{\rL^2(\Omega)}^2  
	\end{equation*}
	and as for $\IV$ we obtain the estimate
	\begin{equation*}
		\int_0^t |\VI|
		\leq \| \partial_z v_{\rE} \|_{\rL^\infty(-h,0)} 
		\cdot \int_0^t \| \nabla v_d \|_{\rL^2(\Omega)}^2
		\leq C \cdot \| \partial_z v_{\rE} \|_{\rL^\infty(-h,0)} \cdot \| v_0 - v_{\rE} \|_{\rL^2(\Omega)}^2 .
	\end{equation*}
	Using $\rH^1(-h,0) \hookrightarrow \rL^\infty(-h,0)$ and the Poincare type inequality we have
	\begin{align*}
		\| w_d \|_{\rL^\infty_z((-h,0),\rL^2_\rH(\T^2))}
		\leq \| w_d \|_{\rH^1_z((-h,0),\rL^2_\rH(\T^2))}
		&\leq C \cdot \| \partial_z w_d \|_{\rL^2_z((-h,0),\rL^2_\rH(\T^2))} \\
		&= C \cdot \| \divH v_d \|_{\rL^2(\Omega)}
		\leq C \cdot \| v_d \|_{\rH^1(\Omega)} 
	\end{align*}
	and conclude for the second term
	\begin{align*}
		|\VII|
		&\leq \| w_d \|_{\rL^\infty_z((-h,0),\rL^2_\rH(\T^2))}
		\cdot \| \partial_z^2 v_{\rE} \|_{\rL^2(-h,0)}
		\cdot \| \partial_z v_d \|_{\rL^2(\Omega)}^2 \\
		&\leq C \cdot \| \partial_z^2 v_{\rE} \|_{\rL^2(-h,0)} \cdot \| v_d \|_{\rH^1(\Omega)}^2 .
	\end{align*}
	Integrating in time yields
	\begin{equation*}
		\int_0^t |\VII|
		\leq C \cdot \| \partial_z^2 v_{\rE} \|_{\rL^2(-h,0)} \cdot \int_0^t \| v_d \|_{\rH^1(\Omega)}^2 
		\leq  C \cdot \| \partial_z^2 v_{\rE} \|_{\rL^2(-h,0)} \cdot \| v_0 - v_{\rE} \|_{\rL^2(\Omega)}^2 .
	\end{equation*}
	Combining the estimates we obtain
	\begin{equation*}
		\begin{aligned}
			&\frac{1}{2} \cdot \| \partial_z v_d(t) \|_{\rL^2(\Omega)}^2 
			+ \nu_{\rH} \int_0^t \| \nablaH \partial_z v_d \|_{\rL^2(\Omega)}^2+ \nu_z \int_0^t \| \partial_z^2 v_d \|_{\rL^2(\Omega)}
			\\ 
			\leq \ &\frac{1}{2} \cdot \| \partial_z (v_0 - v_{\rE}) \|_{\rL^2(\Omega)}^2
			+ 
			C \cdot (1 + \| \partial_z v_{\rE} \|_{\rL^\infty(-h,0)} \\
			+ \ &\| \partial_z^2 v_{\rE} \|_{\rL^2(-h,0)}) \cdot \| v_0 - v_{\rE} \|_{\rL^2(\Omega)}^2 + \varepsilon_1 \cdot \int_0^t \| \nablaH \pi \|_{\rL^2(\T^2)}^2 \\
			+ \ &\left( \varepsilon_3 + \frac{\varepsilon_2}{4 \rho_0 \cdot \varepsilon_1} \right) \cdot \int_0^t  \| \partial_z v_d\|_{\rH^{1}(\Omega)}^2
			+ C \cdot \int_0^t \| \nabla v_d\|_{\rL^2(\Omega)}^2
			\\
			+ \ &C \cdot \int_0^t (\| \nablaH \vbar_d \|_{\rL^2(\Omega)} + \| \nablaH \vbar_d \|_{\rL^2(\Omega)}^2) \cdot \| \partial_z v_d \|_{\rL^2(\Omega)}^2 \, \d s 
			\\
			+ \ &C \cdot \int_0^t \| |\vtilde_d| |\nabla\vtilde_d| \|_{\rL^2(\Omega)}^2 + M_z \cdot \| |\vtilde_d| |\nabla\vtilde_d| \|_{\rL^2(\Omega)}^2.
		\end{aligned}
	\end{equation*}
	Choosing $\varepsilon_3 + \frac{\varepsilon_2}{4 \rho_0 \cdot \varepsilon_1} \leq \frac{1}{2} \min \{\nu_{\rH},\nu_z\}$ we obtain by absorbing 
	\begin{equation*}
		\begin{aligned}
			&\| \partial_z v_d(t) \|_{\rL^2(\Omega)}^2 
			+ \nu_{\rH} \int_0^t \| \nablaH \partial_z v_d \|_{\rL^2(\Omega)}^2
			+ \nu_z \int_0^t \| \partial_z^2 v_d \|_{\rL^2(\Omega)} \\ 
			\leq \ &\| \partial_z (v_0 - v_{\rE}) \|_{\rL^2(\Omega)}^2
			+ 2 \varepsilon_1 \cdot \int_0^t \| \nablaH \pi \|_{\rL^2(\T^2)}^2
			+ C \cdot \int_0^t \| \nabla v_d\|_{\rL^2(\Omega)}^2
			\\
			+ 
			\ &C \cdot (1 + \| \partial_z v_{\rE} \|_{\rL^\infty(-h,0)} + \| \partial_z^2 v_{\rE} \|_{\rL^2(-h,0)}) \cdot \| v_0 - v_{\rE} \|_{\rL^2(\Omega)}^2
			\\
			+ \ &M_z \cdot \| |\vtilde_d| |\nabla\vtilde_d| \|_{\rL^2(\Omega)}^2
			+ C \cdot \int_0^t (\| \nablaH \vbar_d \|_{\rL^2(\Omega)} + \| \nablaH \vbar_d \|_{\rL^2(\Omega)}^2) \cdot \| \partial_z v \|_{\rL^2(\Omega)}^2 \, \d s \\
			+ \ &C \cdot \int_0^t \| |\vtilde_d| |\nabla\vtilde_d| \|_{\rL^2(\Omega)}^2 
		\end{aligned}
	\end{equation*}
	which proves the claim.
\end{proof}
Now, we split \eqref{eq:difference_2} into two evolution equations: one for $\vbar$ and one for $\vtilde$. 
Using
\begin{align*}
	\overline{w_d \cdot  \dz v_{\rE}}
	&= \overline{w_d \cdot \dz \vtilde_E}
	= - \overline{\dz w_d \cdot \vtilde_E}
	= \overline{\divH v_d \cdot \vtilde_E}  
	= \overline{\divH \vtilde_d \cdot \vtilde_E}  , \\
	\overline{w_d \cdot  \dz v_d}
	&= \overline{w_d \cdot \dz \vtilde_d}
	= - \overline{\dz w_d \cdot \vtilde_d}
	= \overline{\divH v_d \cdot \vtilde_d}  
	= \overline{\divH \vtilde_d \cdot \vtilde_d}  
\end{align*}
we obtain

\begin{equation}
	\label{eq:difference_2 bar}
	\left\{
	\begin{aligned}
		\partial_t \vbar_d - \nu_{\rH} \DeltaH \vbar_d + \frac{1}{\rho_0}\nablaH \pi_d + f \vbar_d^\perp
		&=
		\ - \vbar_d \cdot \nabla \vbar_d - \vbar_{\rE} \cdot \nablaH \vbar_d
		-\frac{1}{h} \nu_z \partial_z v_d\big|_{\Gamma_b} \\ 
		&\phantom{ = } \ \ - \overline{\vtilde_d \cdot \nablaH \vtilde_d + \divH \vtilde_d \cdot \vtilde_d} \\ 
		&\phantom{ = } \ \ 
		- \overline{\vtilde_{\rE} \cdot \nablaH \vtilde_d + \divH \tilde{v}_d \cdot v_{\rE}}, 
		\\
		\divH \vbar_d &= 0, 
	\end{aligned}
	\right.
\end{equation}
in $(0,T) \times \T^2$ and 
\begin{equation}
	\label{eq:difference_2 tilde} 
	\begin{aligned}
		&\phantom{ = } \ \ \partial_t \vtilde_d - \nu_{\rH} \DeltaH \vtilde_d - \nu_z \partial_z^2 \vtilde_d + \vtilde_d \cdot \nablaH \vtilde_d + w_d \partial_z \vtilde_d +\vbar_d \cdot \nablaH \vtilde_d + f \vtilde_d^\perp \\
		&= - \vtilde_d \cdot \nablaH \vbar_d
		+ \overline{\vtilde_d \cdot \nablaH \vtilde_d + \divH \tilde{v}_d \cdot v_d}
		+ \frac{1}{h} \nu_z \partial_z v_d\big|_{\Gamma_b} - w_d \partial_z v_{\rE} \\
		&\phantom{ = } \
		- \vtilde_{\rE} \cdot \nablaH \vtilde_d
		- \vtilde_{\rE} \cdot \nablaH \vbar_d
		- \vbar_{\rE} \cdot \nablaH \vtilde_d 
		+ \overline{\vtilde_{\rE} \cdot \nablaH \vtilde_d + \divH \tilde{v}_d \cdot v_{\rE}}, \\
	\end{aligned}
\end{equation}
in $(0,T) \times \Omega$
subject to the boundary conditions
\begin{equation*}
	\partial_z \vtilde_d = 0 \enspace \text{ on } \enspace (0,T) \times \Gamma_u, \enspace \text{ and } \enspace
	\vtilde_d = -\vbar \enspace \text{ on } \enspace (0,T) \times \Gamma_b.
\end{equation*}

Note that the left-hand side of \eqref{eq:difference_2 bar} is a two-dimensional Stokes equation. As in \cite[Lemma 6.1]{HK:16_2} we use the following well-known estimate of the Stokes operator in order to establish bounds for \eqref{eq:difference_2 bar}.

\begin{lem}[Estimates of the 2D Stokes equations]\label{lem:2D NSE_2}
	Let $(u,\pi)$ be the solution of the two-dimensional Stokes equations with rotation
	\begin{equation*} 
		\left\{
		\begin{aligned}
			\partial_t u - \nu_{\rH} \DeltaH u + \nablaH \pi + u^\perp &= k,\\
			\divH u &= 0
		\end{aligned}
		\right. 
	\end{equation*} 
	for $k \in \rL^2(0,T;\rL^2(\T^2))$. Then $(u,\pi)$ satisfies the estimate 
	\begin{equation*}
		\sqrt{\nu_{\rH}} \cdot \partial_t \| \nablaH u \|_{\rL^2(\T^2)}^2
		+ \nu_{\rH} \cdot \| \DeltaH u \|_{\rL^2(\T^2)}^2
		+ \| \nablaH \pi \|_{\rL^2(\T^2)}^2
		\leq C \cdot (\| h \|_{\rL^2(\T^2)}^2+\| u \|_{\rL^2(\T^2)}^2) .
	\end{equation*}
\end{lem}
\begin{proof}
	Multiplying by $\partial_t u - \DeltaH u$, integration by parts and Young's inequality yield
	\begin{align*}
		\| \partial_t u \|_{\rL^2(\T^2)}^2 
		+
		\sqrt{\nu_{\rH}}\partial_t \| \nablaH u \|_{\rL^2(\T^2)}^2
		+\nu_{\rH} \|\DeltaH u \|_{\rL^2(\T^2)}^2
		\leq \int_{\T^2} k (\partial_t u - \DeltaH u) \\
		\leq 
		C \| k \|_{\rL^2(\T^2)}^2
		+
		\varepsilon \| \partial_t u \|_{\rL^2(\T^2)}^2 
		+
		\varepsilon 
		\|\DeltaH u \|_{\rL^2(\T^2)}^2
		, 
	\end{align*}
	where we have used 
	\begin{align*}
		-\int_{\T^2} \nablaH \pi \DeltaH u
		&=  \int_{\T^2} \pi \divH \DeltaH u
		= \int_{\T^2} \pi \DeltaH \divH u
		= 0, \\
		-\int_{\T^2} u \DeltaH u 
		&= \int_{\T^2} \nablaH  u^\perp \nablaH u
		= \int_{\T^2} (\nablaH  u)^\perp \nablaH u
		= 0 . 
	\end{align*}
	By absorbing this implies 
	\begin{equation*}
		\| \partial_t u \|_{\rL^2(\T^2)}^2 
		+
		\sqrt{\nu_{\rH}}\cdot \partial_t \| \nablaH u \|_{\rL^2(\T^2)}^2
		+ \nu_{\rH}\cdot \|\DeltaH u \|_{\rL^2(\T^2)}^2
		\leq 
		C \cdot \| k \|_{\rL^2(\T^2)}^2
		, 
	\end{equation*}
	The pressure satisfies the elliptic equation
	\begin{equation*}
		\DeltaH \pi = \divH k - \divH u^\perp  
	\end{equation*}
	which yields to the estimate
	\begin{align*}
		\| \nablaH \pi \|_{\rL^2(\Omega)}
		&\leq C \cdot (\| \divH k \|_{\rH^{-1}(\T^2)}
		+ \| \divH u^\perp \|_{\rH^{-1}(\T^2)}) \\
		&\leq C \cdot (\| k \|_{\rL^2(\T^2)}
		+ \| u \|_{\rL^2(\T^2)}) ,
	\end{align*}
	where we have used the characterization of the negative order Sobolev space
	\begin{equation*}
		\rH^{-1}(\T^2) = \{ \divH k \colon k \in \rL^2(\T^2) \} .
	\end{equation*}
	Adding both estimates yields the claim. 
\end{proof}

Now we establish estimates for the barotropic modes $\vbar_d$.

\begin{lem}[Estimates for $\bar{v}_d$]\label{lem:estimate vbar_2}
	Suppose $\alpha < 0$. 
	Then 
	\begin{equation*}
		\begin{aligned}
			&\sqrt{\nu_{\rH}}\cdot \| \nablaH \vbar_d(t) \|_{\rL^2(\T^2)}^2
			+ \nu_{\rH} \cdot \int_0^t \| \DeltaH \vbar_d \|_{\rL^2(\T^2)}^2
			+ \int_0^t \| \nablaH \pi \|_{\rL^2(\T^2)}^2 \\ 
			\leq \ &
			C \cdot (1 + \| v_{\rE} \|_{\rL^\infty(-h,0)}^2) \cdot \| v_0 - v_{\rE} \|_{\rL^2(\Omega)}^2 \\
			+ \ &
			\bar{\varepsilon} \cdot \int_0^t \| \nabla \partial_z v \|_{\rL^2(\Omega)}^2
			+ \bar{M} \cdot \int_0^t \| |\vtilde_d| |\nabla\vtilde_d| \|_{\rL^2(\Omega)}^2 \\
			+ \ &C \cdot \int_0^t \bigl( 1+ \| v_d \|_{\rL^2(\Omega)}^2 + \| v_d \|_{\rL^2(\Omega)}^4 \bigr) \cdot \| v_d \|_{\rH^1(\Omega)}^2 \\
			+ \ &C \cdot \int_0^t \bigl( 1+ \| v_d \|_{\rL^2(\Omega)} + \| v_d \|_{\rL^2(\Omega)}^2 \bigr) \cdot \bigl( \| v_d \|_{\rH^1(\Omega)}+\| v_d \|_{\rH^1(\Omega)}^2 \bigr) \cdot \| \nablaH \vbar_d \|_{\rL^2(\T^2)}^2 
		\end{aligned}
	\end{equation*}
	for $\bar{\varepsilon} > 0$ and a fixed constant $\bar{M} > 0$.
\end{lem}
\begin{proof}
	As in \cite[Step 1]{HK:16_2} or \cite[Step 2]{GGHHK:20_2} we 
	apply \autoref{lem:2D NSE_2} to \eqref{eq:difference_2 bar} and obtain
	\begin{equation*}
		\begin{aligned}
			&\sqrt{\nu_{\rH}}\cdot \partial_t \| \nablaH \vbar_d \|_{\rL^2(\T^2)}^2
			+ \nu_{\rH} \cdot \| \DeltaH \vbar_d \|_{\rL^2(\T^2)}^2
			+ \| \nablaH \pi \|_{\rL^2(\T^2)}^2
			\\
			\leq \ &C \cdot \bigl(\| \vbar_d \nablaH \vbar_d \|_{\rL^2(\T^2)}^2
			+ \| |\vtilde_d| |\nablaH \vtilde_d| \|_{\rL^2(\Omega)}^2
			+ \| \partial_z v_d\big|_{\Gamma_b} \|_{\rL^2(\T^2)}^2\bigr) \\
			+ \ &C \cdot
			\bigl( 
			\| \vbar_{\rE} \nablaH \vbar_d \|_{\rL^2(\T^2)}^2
			+ \| |\vtilde_{\rE}| |\nablaH \vtilde_d| \|_{\rL^2(\Omega)}^2
			\bigr) 
			+ C \cdot \| \vbar_d \|_{\rL^2(\T^2)}^2\\
			=: \ & \ \I + \II + \III + \IV + \V + \VI .
		\end{aligned}
	\end{equation*}
	From \cite[Step 1, (6.5)]{HK:16_2} we have estimates for the first three terms
	\begin{align*}
		\I + \II + \III 
		&\leq
		\varepsilon \cdot \| \nabla \partial_z v_d \|_{\rL^2(\Omega)}^2
		+ M_z \cdot \| |\vtilde_d| |\nabla \vtilde_d| \|_{\rL^2(\Omega)}^2\\
		&+ C \cdot \bigl( 1+ \| v_d \|_{\rL^2(\Omega)}^2 + \| v_d \|_{\rL^2(\Omega)}^4 \bigr) \cdot \| v_d \|_{\rH^1(\Omega)}^2 \\
		&+ C \cdot \bigl( 1+ \| v_d \|_{\rL^2(\Omega)} + \| v_d \|_{\rL^2(\Omega)}^2 \bigr) 
		\\ &\hspace{24pt} \cdot \bigl( \| v_d \|_{\rH^1(\Omega)}+\| v_d \|_{\rH^1(\Omega)}^2 \bigr) \cdot \| \nablaH \vbar_d \|_{\rL^2(\T^2)}^2
	\end{align*}
	for $\varepsilon > 0$. 
	For the last term we obtain by \autoref{lem:energy_2} that
	\begin{equation*}
		\int_0^t \VI
		\leq \int_0^t \| v_d \|_{\rL^2(\Omega)}^2
		\leq C \cdot \| v_0 - v_{\rE} \|_{\rL^2(\Omega)}^2 .
	\end{equation*}
	For the two remaining terms we obtain
	\begin{equation*}
		\IV + \V 
		\leq 2 
		\| v_{\rE} \nablaH v_d \|_{\rL^2(\Omega)}^2
		\leq 2 
		\| v_{\rE} \|_{\rL^\infty(-h,0)}^2
		\cdot \| \nablaH v_d \|_{\rL^2(\Omega)}^2 .
	\end{equation*}
	Integrating in time and using \autoref{lem:energy_2} yields
	\begin{equation*}
		\int_0^t \IV + \V 
		\leq 2 
		\| v_{\rE} \|_{\rL^\infty(-h,0)}^2
		\cdot \int_0^t \| \nablaH v_d \|_{\rL^2(\Omega)}^2 
		\leq C 
		\cdot
		\| v_{\rE} \|_{\rL^\infty(-h,0)}^2 
		\cdot 
		\| v_0 - v_{\rE} \|_{\rL^2(\Omega)}^2 .
		\qedhere
	\end{equation*}
\end{proof}

Finally, we show estimates for the baroclinic modes $\vtilde_d$.

\begin{lem}[Estimates for $\vtilde_d$]\label{lem:estimate vtilde_2}
	Suppose $\alpha < 0$. 
	Then 
	\begin{equation*}
		\begin{aligned}
			&\| \vtilde_d(t) \|_{\rL^4(\Omega)}^4 
			+ \int_0^t \| \nabla |\vtilde_d|^2 \|_{\rL^2(\Omega)}^2
			+ \int_0^t \| |\vtilde_d| |\nabla \vtilde_d| \|_{\rL^2(\Omega)}^2 \\
			\leq \, & C \cdot \| v_0-v_{\rE} \|_{\rL^4(\Omega)}^4
			+ C \cdot 
			\| \dz v_{\rE} \|_{\rL^\infty(-h,0)}^2 \cdot \| v_0 - v_{\rE} \|_{\rL^2(\Omega)}^2 \\
			&+ C \cdot \int_0^t \bigl( \| v_d \|_{\rH^1(\Omega)}^{\frac{2}{3}} + \| v_d \|_{\rH^1(\Omega)}+\| v_d \|_{\rH^1(\Omega)}^2 \bigr) 
			\cdot \| \vtilde_d \|_{\rL^4(\Omega)}^4  \\
			&+ \frac{\tilde{\varepsilon}}{\tilde{M}} \cdot \int_0^t \| \nablaH \partial_z v \|_{\rL^2(\Omega)}^2 . 
		\end{aligned}
	\end{equation*}
	for $\tilde{\varepsilon} > 0$ and $\tilde{M} := 2(\bar{M}+2 M_z) > 0$.
\end{lem}
\begin{proof}
	Following \cite[Step 3]{HK:16_2} or \cite[Step 1]{GGHHK:20_2} we multiply by $|\vtilde_d|^2 \vtilde_d$ and integrate 
	\begin{equation*}
		\begin{aligned}
			&  \frac{1}{4} \, \partial_t \| \vtilde_d \|_{\rL^4(\Omega)}^4 
			+ \frac{1}{2} \, \| \nabla |\vtilde_d|^2 \|_{\rL^2(\Omega)}^2
			+ \| |\vtilde_d| |\nabla \vtilde_d| \|_{\rL^2(\Omega)}^2
			\\ 
			= \ &-\int_{\Omega} (\vtilde_d \cdot \nablaH \vbar_d) \cdot |\vtilde_d|^2 \vtilde_d + \frac{1}{h} \int_{\Omega} \int_{-h}^0 (\vtilde_d \cdot \nablaH \vtilde_d + \divH \vtilde_d \cdot \vtilde_d) \d z \cdot |\vtilde_d|^2 \vtilde_d  
			\\
			\ &+ \frac{1}{h} \int_{\Omega} (\partial_z v_d)\big|_{\Gamma_b} \cdot |\vtilde_d|^2 \vtilde_d
			-
			\int_{\Omega} w_d \partial_z v_{\rE} \cdot |\vtilde_d|^2 \vtilde_d 
			-
			\int_{\Omega} \vtilde_{\rE} \cdot \nablaH \vbar_d \cdot |\vtilde_d|^2 \vtilde_d \\
			&- \frac{1}{h} \int_{\Omega} \int_{-h}^0 (\vtilde_{\rE} \cdot \nablaH \vtilde_d + \divH \vtilde_d \cdot \vtilde_{\rE}) \d z \cdot |\vtilde_d|^2 \vtilde_d \\
			=: \ & \I + \II + \III + \IV + \V + \IV , 
		\end{aligned}
	\end{equation*}
	where we have used the following cancellation laws (see \cite[Lemma 6.3]{HK:16_2})
	\begin{equation*} 
		\begin{aligned}
			&\int_{\Omega} \vtilde_d \cdot \nablaH \vtilde_d \cdot |\vtilde_d|^2 \vtilde_d  + \int_{\Omega} w_d \partial_z \vtilde_d \cdot |\vtilde_d|^2 \vtilde_d = 0, \hspace{-0.5em} 
			&&\int_{\Omega} \vbar_d \cdot \nablaH \vtilde_d \cdot |\vtilde_d|^2 \vtilde_d = 0, \\
			&\int_{\Omega} \vtilde_{\rE} \cdot \nablaH \vtilde_d \cdot |\vtilde_d|^2 \vtilde_d = 0, &&\int_{\Omega} \vbar_{\rE} \cdot \nablaH \vtilde_d \cdot |\vtilde_d|^2 \vtilde_d = 0 ,
		\end{aligned}
	\end{equation*}
	and the identity
	\begin{equation*}
		\int_{\Omega} f \vtilde_d^\perp \cdot |\vtilde_d|^2 \vtilde_d = 0 .
	\end{equation*}
	From \cite[Step 3, (6.9)]{HK:16_2} we know that
	\begin{align*}
		\I + \II + \III
		&\leq C \cdot  \bigl( \| v_d \|_{\rH^1(\Omega)}^{\frac{2}{3}} + \| v_d \|_{\rH^1(\Omega)}+\| v_d \|_{\rH^1(\Omega)}^2 \bigr) 
		\cdot \| \vtilde_d \|_{\rL^4(\Omega)}^4   \\
		&+ \frac{\tilde{\varepsilon}_1}{\tilde{M}} \cdot \| \nablaH \partial_z v \|_{\rL^2(\Omega)}^2 
		+ \tilde{\varepsilon}_2 \cdot \| \nabla |\vtilde_d|^2 \|_{\rL^2(\Omega)}^2
	\end{align*}
	for $\tilde{\varepsilon}_1,\tilde{\varepsilon}_2 > 0$ and $\tilde{M} = 2(\bar{M} + 2 M_z) > 0$.
	Using the embeddings
	\begin{align*}
		\rH^1(-h,0) \hookrightarrow \rL^\infty(-h,0), \qquad 
		\rH^1(\T^2) \hookrightarrow \rL^4(\T^2), \\
		\rL^4(\Omega) \hookrightarrow \rL^3_z((-h,0); \rL^4(\T^2)) \hookrightarrow  \rL^2_z((-h,0); \rL^2(\T^2)),
	\end{align*}
	we estimate the forth term by 
	\begin{align*}
		|\IV|
		&\leq \| \dz v_{\rE} \|_{\rL^\infty(-h,0)} \int_{\T^2} \| w_d \|_{\rL^\infty_z(-h,0)} \int_{-h}^0 |\vtilde_d|^3 
		\\
		&\leq \| \dz v_{\rE} \|_{\rL^\infty(-h,0)} \cdot \| \divH v_d \|_{\rL^2(\Omega)} \cdot  \int_{-h}^0 \| |\vtilde_d|^3 \|_{\rL^2(\T^2)} \\
		&\leq \| \dz v_{\rE} \|_{\rL^\infty(-h,0)} \cdot \| v_d \|_{\rH^1(\Omega)} \cdot \int_{-h}^0 \bigl( \| \vtilde_d \|_{\rL^4(\T^2)}^3 + \| \vtilde_d \|_{\rL^4(\T^2)} \cdot \| \nablaH |\vtilde_d|^2 \|_{\rL^2(\T^2)} \bigr) \\
		&\leq \| \dz v_{\rE} \|_{\rL^\infty(-h,0)} \cdot \| v_d \|_{\rH^1(\Omega)} \cdot \| \vtilde_d \|_{\rL^4(\Omega)}^3 \\ 
		&+ \| \dz v_{\rE} \|_{\rL^\infty(-h,0)} \cdot \| v_d \|_{\rH^1(\Omega)} \cdot \| \vtilde_d \|_{\rL^4(\Omega)} \cdot \| \nablaH |\vtilde_d|^2 \|_{\rL^2(\Omega)} \\
		&\leq
		\| \dz v_{\rE} \|_{\rL^\infty(-h,0)} \cdot \| v_d \|_{\rH^1(\Omega)}
		+ \| \dz v_{\rE} \|_{\rL^\infty(-h,0)} \cdot \| v_d \|_{\rH^1(\Omega)} \cdot \| \vtilde_d \|_{\rL^4(\Omega)}^4 \\
		&\hspace{10pt} + C \cdot \| \dz v_{\rE} \|_{\rL^\infty(-h,0)}^2 \cdot \| v_d \|_{\rH^1(\Omega)}^2 \cdot \| \vtilde_d \|_{\rL^4(\Omega)}^2 + \tilde{\varepsilon}_3 \cdot 
		\| \nablaH |\vtilde_d|^2 \|_{\rL^2(\Omega)}^2 \\
		&\leq 1 + C \cdot 
		\| \dz v_{\rE} \|_{\rL^\infty(-h,0)}^2 \cdot \| v_d \|_{\rH^1(\Omega)}^2 \\
		&\hspace{10pt} + \| \dz v_{\rE} \|_{\rL^\infty(-h,0)} \cdot \| v_d \|_{\rH^1(\Omega)} \cdot \| \vtilde_d \|_{\rL^4(\Omega)}^4 
		+ \tilde{\varepsilon}_3 \cdot 
		\| \nablaH |\vtilde_d|^2 \|_{\rL^2(\Omega)}^2 .
	\end{align*}
	We argue in a similar way for the next term
	\begin{align*}
		| \V | &\leq \| \vtilde_{\rE}\|_{\rL^\infty(-h,0)} \cdot \int_{\T^2} |\nablaH \vbar_d | \cdot \int_{-h}^0 |\vtilde_d|^3 \\
		&\leq \| \vtilde_{\rE}\|_{\rL^\infty(-h,0)} \cdot  \| v_d \|_{\rH^1(\Omega)} \cdot \int_{-h}^0  \| |\vtilde_d|^3 \|_{\rL^2(\T^2)}
	\end{align*}
	and now it follows exactly as above that
	\begin{align*}
		| \V |
		\leq 1 + C \cdot 
		\| \dz v_{\rE} \|_{\rL^\infty(-h,0)}^2 \cdot \| v_d \|_{\rH^1(\Omega)}^2 + \| \dz v_{\rE} \|_{\rL^\infty(-h,0)} \cdot \| v_d \|_{\rH^1(\Omega)} \cdot \| \vtilde_d \|_{\rL^4(\Omega)}^4 \\
		+ \tilde{\varepsilon}_4 \cdot 
		\| \nablaH |\vtilde_d|^2 \|_{\rL^2(\Omega)}^2 .
	\end{align*}
	Finally, we get for the last term
	\begin{align*}
		|\VI|
		&\leq \frac{1}{h} \int_{\T^2} \int_{-h}^0 |\vtilde_{\rE} \cdot \nablaH \vtilde_d + \divH \vtilde_d \cdot \vtilde_{\rE}| \cdot \int_{-h}^0 |\vtilde_d|^3 \\
		&\leq \frac{1}{h} \int_{-h}^0 \| \vtilde_{\rE} \cdot \nablaH \vtilde_d + \divH \vtilde_d \cdot \vtilde_{\rE} \|_{\rL^2(\T^2)} \cdot \int_{-h}^0 \| |\vtilde_d|^3 \|_{\rL^2(\T^2)} \\
		&\leq \frac{1}{h} \cdot \| \vtilde_{\rE} \cdot \nablaH \vtilde_d + \divH \vtilde_d \cdot \vtilde_{\rE} \|_{\rL^2(\Omega)} \cdot \int_{-h}^0 \| |\vtilde_d|^3 \|_{\rL^2(\T^2)} \\
		&\leq \frac{2}{h} \cdot \| \vtilde_{\rE} \|_{\rL^\infty(-h,0)} \cdot \| \vtilde_d \|_{\rH^1(\Omega)} \cdot \int_{-h}^0 \| |\vtilde_d|^3 \|_{\rL^2(\T^2)}
	\end{align*}
	and we conclude as before that
	\begin{align*}
		| \VI |
		\leq 1 &+ C \cdot 
		\| \dz v_{\rE} \|_{\rL^\infty(-h,0)}^2 \cdot \| v_d \|_{\rH^1(\Omega)}^2 \\
		&+ \| \dz v_{\rE} \|_{\rL^\infty(-h,0)} \cdot \| v_d \|_{\rH^1(\Omega)} \cdot \| \vtilde_d \|_{\rL^4(\Omega)}^4
		+ \tilde{\varepsilon}_5 \cdot 
		\| \nablaH |\vtilde_d|^2 \|_{\rL^2(\Omega)}^2 .
	\end{align*}
	Using \autoref{lem:energy_2} we see that the first term is bounded by
	\begin{equation*}
		C \cdot
		\| \dz v_{\rE} \|_{\rL^\infty(-h,0)}^2 \cdot \int_0^t \| v_d \|_{\rH^1(\Omega)}^2
		\leq C \cdot 
		\| \dz v_{\rE} \|_{\rL^\infty(-h,0)}^2 \cdot \| v_0 - v_{\rE} \|_{\rL^2(\Omega)}^2 .
	\end{equation*}
	Combining all estimates, integrating in time and absorbing the term \\ $\| \nablaH |\vtilde_d|^2 \|_{\rL^2(\Omega)}^2$ we obtain
	\begin{align*}
		&\| \vtilde_d(t) \|_{\rL^4(\Omega)}^4 
		+ \int_0^t \| \nabla |\vtilde_d|^2 \|_{\rL^2(\Omega)}^2
		+ \int_0^t \| |\vtilde_d| |\nabla \vtilde_d| \|_{\rL^2(\Omega)}^2 \\
		\leq \, & C \cdot \| v_0-v_{\rE} \|_{\rL^4(\Omega)}^4
		+ C \cdot 
		\| \dz v_{\rE} \|_{\rL^\infty(-h,0)}^2 \cdot \| v_0 - v_{\rE} \|_{\rL^2(\Omega)}^2 \\
		&+ C \cdot \int_0^t \bigl( \| v_d \|_{\rH^1(\Omega)}^{\frac{2}{3}} + \| v_d \|_{\rH^1(\Omega)}+\| v_d \|_{\rH^1(\Omega)}^2 \bigr) 
		\cdot \| \vtilde_d \|_{\rL^4(\Omega)}^4  \\
		&+ \frac{\tilde{\varepsilon}}{\tilde{M}} \cdot \int_0^t \| \nablaH \partial_z v_d \|_{\rL^2(\Omega)}^2 . \qedhere 
	\end{align*} 
\end{proof}

Combining our previous estimates we obtain the following result. 

\begin{prop}[$\rL^\infty_t\rH^1$-$\rL^2_t \rH^2$-estimates]\label{prop:estimate_2}
	Suppose $\alpha < 0$. 
	Let $v_d$ be a solution to \eqref{eq:difference_2}. Then there exists a bound $B$ continuously depending on $t, \| v_0 \|_{\rH^1(\Omega)}$ and $\| v_{\rE} \|_{\rH^2(\Omega)}$ such that 
	\begin{equation*}
		\| v_d(t) \|_{\rH^1(\Omega)}^2
		+ \int_0^t \| v_d(s) \|_{\rH^2(\Omega)}^2 \d s \leq B(t) .
	\end{equation*}
\end{prop}
\begin{proof}
	We follow the strategy from \cite[Step 4]{HK:16_2} or \cite[Step 4]{GGHHK:20_2}. Adding the estimates from \autoref{lem:estimate vz_2}, \autoref{lem:estimate vbar_2} and \autoref{lem:estimate vtilde_2} and absorbing the terms $\| |\vtilde_d| |\nabla \vtilde_d| \|_{\rL^2(\Omega)}^2$, $\| \nabla \partial_z v_d \|_{\rL^2(\Omega)}^2$ and $\| \nablaH \pi \|_{\rL^2(\Omega)}^2$ we obtain
	\begin{align*}
		&\| \nablaH \vbar_d(t) \|_{\rL^2(\T^2)}^2
		+ \| \partial_z v_d(t) \|_{\rL^2(\Omega)}^2
		+ \tilde{M} \cdot \| \vtilde_d(t) \|_{\rL^4(\Omega)}^4
		\\
		\enspace &
		+\int_0^t \| \nablaH \pi\|_{\rL^2(\T^2)}^2
		+ \| \nabla \partial_z v_d \|_{\rL^2(\Omega)}^2
		+ \tilde{M} \cdot \| |\vtilde_d| |\nablaH \vtilde_d| \|_{\rL^2(\Omega)}^2 
		\\
		\leq \enspace &\int_0^t K_1(s) \cdot \bigl( \| \nablaH \vbar_d(t) \|_{\rL^2(\T^2)}^2
		+ \| \partial_z v(t) \|_{\rL^2(\Omega)}^2
		+ \tilde{M} \cdot \| \vtilde_d(t) \|_{\rL^4(\Omega)}^4 \bigr) + K_2(s)
		+ K_0 ,
	\end{align*}
	where 
	\begin{align*}
		K_0 &:= C \| \partial_z (v_0 - v_{\rE}) \|_{\rL^2(\Omega)}^2
		+ C (1 + \| \partial_z v_{\rE} \|_{\rL^\infty(-h,0)} \\
		&+ \| \partial_z v_{\rE} \|_{\rL^\infty(-h,0)}^2 + \| \partial_z^2 v_{\rE} \|_{\rL^2(-h,0)}) \cdot \| v_0 - v_{\rE} \|_{\rL^2(\Omega)}^2 \\
		K_1(t) &:= C \cdot \bigl(1+ \| v_d(t) \|_{\rL^2(\Omega)} +\| v_d(t) \|_{\rL^2(\Omega)}^2\bigr) \\
		&\hspace{27pt}\cdot \bigl( \| v_d(t) \|_{\rH^1(\Omega)}^{\frac{2}{3}} + \| v_d (t) \|_{\rH^1(\Omega)} + \| v_d(t) \|_{\rH^1(\Omega)}^2 \bigr),\\
		K_2(t) &:= C \cdot \bigl(1+ \| v_d(t) \|_{\rL^2(\Omega)}^2 +\| v_d(t) \|_{\rL^2(\Omega)}^4\bigr) \cdot \| v_d(t) \|_{\rH^1(\Omega)}^2 .
	\end{align*}
	Now Hölder's inequality and the energy estimate from \autoref{lem:energy_2} imply $K_1, K_2 \in \rL^1(0,T)$ and Gronwall's inequality implies the estimate
	\begin{equation}
		\label{eq:bound_2}
		\begin{aligned}
			&\| \nablaH \vbar_d(t) \|_{\rL^2(\T^2)}^2
			+ \| \partial_z v(t) \|_{\rL^2(\Omega)}^2
			+ \| \vtilde_d(t) \|_{\rL^4(\Omega)}^4
			\\
			+ \
			&\int_0^t \| \nablaH \pi\|_{\rL^2(\T^2)}^2
			+ \| \nabla \partial_z v \|_{\rL^2(\Omega)}^2
			+ \| |\vtilde_d| |\nablaH \vtilde_d| \|_{\rL^2(\Omega)}^2 
			\leq \hat{B},  
		\end{aligned}
	\end{equation}
	where the bound $\hat{B} := \hat{B}(t,\| v_0 \|_{\rH^1(\Omega)},\|v_{\rE}\|_{\rH^2(-h,0)} )$ is given by
	\begin{equation*}
		\hat{B}(t,\| v_0 	\|_{\rH^1(\Omega)},\|v_{\rE}\|_{\rH^2(-h,0)} ) := 2 \cdot \max \{ 1, \nicefrac{1}{\tilde{M}} \} \cdot \left( K_0 + \int_0^t K_2(s) \, \d s \right)
		\mathrm{e}^{ \int_0^t K_1(s) \, \d s } .
	\end{equation*}
	We write \eqref{eq:difference_2}$_1$ as
	\begin{equation*}
		\partial_t v_d -\nu_{\rH}\Delta v_d-\nu_z\partial_z^2 v_d
		+ f v_d^\perp = F  
	\end{equation*}
	with
	\begin{equation*}
		F := - u_d \cdot \nabla v_d - v_{\rE} \cdot \nablaH v_d 
		- w_d \dz v_{\rE} - \frac{1}{\rho_0}\nablaH \pi_d.
	\end{equation*}
	Multiplying \eqref{eq:difference_2} by $-\nu_{\rH}\Delta v_d-\nu_z\partial_z v_d$, integrating by parts and Young's inequality yield  
	\begin{align*}
		\nu_{\rH} \partial_t \| \nablaH v_d \|_{\rL^2(\Omega)}^2
		+ \nu_z \partial_t \| \dz v_d \|_{\rL^2(\Omega)}
		+ \| \nu_{\rH}  \DeltaH v_d + \nu_z \partial_z^2 v_d\|_{\rL^2(\Omega)}^2		\\
		\leq C \cdot \| F \|_{\rL^2(\Omega)}^2
		+ \frac{1}{2} \| \nu_{\rH} \DeltaH v_d + \nu_z \partial_z^2 v_d\|_{\rL^2(\Omega)}^2,
	\end{align*}
	where we have used
	\begin{equation*}
		\int_{\Omega} f v_d^\perp (-\nu_{\rH}\Delta v_d-\nu_z\partial_z v_d)
		= - \int_{\Omega} f \nu_{\rH} \nablaH v_d^\perp \nablaH v_d 
		- \int_{\Omega} f \nu_z \partial_z v_d^\perp \partial_z v_d
		= 0 .
	\end{equation*}
	Absorbing the term $\| \nu_{\rH}^2 \DeltaH v_d +\nu_z^2 \partial_z^2 v_d \|_{\rL^2(\Omega)}$ implies
	\begin{equation*}
		\nu_{\rH} \cdot \partial_t \| \nablaH v_d \|_{\rL^2(\Omega)}
		+ \nu_z \partial_t \| \dz v_d \|_{\rL^2(\Omega)}
		+ \|\DeltaH v_d + \nu_z \partial_z^2 v_d\|_{\rL^2(\Omega)}^2
		\leq C \cdot \| F \|_{\rL^2(\Omega)}^2 .
	\end{equation*}
	It remains to bound 
	\begin{align*}
		\| F \|_{\rL^2(\Omega)}^2
		\leq \enspace &C \cdot \bigl(
		\| \vbar_d \cdot \nablaH \vbar_d \|_{\rL^2(\T^2)}^2
		+	\| \vbar_d \cdot \nablaH \vtilde_d \|_{\rL^2(\Omega)}^2 
		+	\| \vtilde_d \cdot \nablaH \vbar_d \|_{\rL^2(\Omega)}^2  \\
		&\hspace{30pt} + \| \vtilde_d \cdot \nablaH \vtilde_d \|_{\rL^2(\Omega)}^2
		+ \| w_d \partial_z v_d \|_{\rL^2(\Omega)}^2
		+ \| \nablaH \pi \|_{\rL^2(\T^2)}^2 \\
		&\hspace{30pt}+ \| v_{\rE} \nabla v_d \|_{\rL^2(\Omega)}^2
		+ \| w_d \partial_z v_{\rE} \|_{\rL^2(\Omega)}^2
		\bigr) \\
		=: \enspace &C \cdot ( \I + \II + \III + \IV + \V + \VI + \VII + \VIII ) .
	\end{align*} 
	Using interpolation estimates and Sobolev embeddings it is shown in \cite[Step~4]{HK:16_2} or \cite[Step 4]{GGHHK:20_2} that
	\begin{align*}
		\I + \II + \III + \IV + \V + \VI \,
		\leq \,  
		&C \cdot \| \partial_z v_d \|_{\rL^2(\Omega)}^2 \cdot \| \nabla v_d \|_{\rL^2(\Omega)}^2 \cdot \| \nabla \partial_z v_d \|_{\rL^2(\Omega)}^2 \\
		&+ C \cdot \bigl( \| \vbar_d \|_{\rH^1(\T^2)}^2 + \| \vtilde_d\|_{\rL^4(\Omega)}^4 \bigr) \cdot \| v_d \|_{\rH^1(\Omega)}^2 .
	\end{align*}
	For the remaining two terms we note that
	\begin{equation*}
		\VII  
		\leq  \| v_{\rE} \|_{\rL^\infty(-h,0)}^2 \cdot \| v_d \|_{\rH^1(\Omega)}^2 
	\end{equation*}
	and
	\begin{equation*}
		\VIII \leq \| \dz v_{\rE} \|_{\rL^\infty(-h,0)}^2 \cdot \| w_d \|_{\rL^2(\Omega)}^2
		\leq \| \dz v_{\rE} \|_{\rL^\infty(-h,0)}^2 \cdot \| v_d \|_{\rH^1(\Omega)}^2 .
	\end{equation*}
	Using the energy estimates, see \autoref{lem:energy_2}, we obtain
	\begin{align*}
		\int_0^t \VII + \VIII 
		&\leq \bigl( \| v_{\rE} \|_{\rL^\infty(-h,0)}^2 + \| \dz v_{\rE} \|_{\rL^\infty(-h,0)}^2 \bigr) \cdot \int_0^t \| v_d \|_{\rH^1(\Omega)}^2 \\
		&\leq C \cdot \bigl(\| v_{\rE} \|_{\rL^\infty(-h,0)}^2 + \| \dz v_{\rE} \|_{\rL^\infty(-h,0)}^2 \bigr) \cdot \| v_{\rE} - v_0 \|_{\rL^2(\Omega)}^2 .
	\end{align*}
	Combining these estimates we obtain
	\begin{align*}
		\| v_d(t) \|_{\rH^1(\Omega)}^2
		+ \int_0^t \| v_d \|_{\rH^2(\Omega)}^2
		&\leq 
		C \cdot \nu_{\rH} \| \nablaH v_0(t) \|_{\rL^2(\Omega)}^2
		+ C \cdot \nu_z \| \partial_z v_d(t) \|_{\rL^2(\Omega)}^2
		\\
		&+ C \int_0^t \| \nu_{\rH}^2 \DeltaH v_d +\nu_z^2 \partial_z^2 v_d \|_{\rL^2(\Omega)} \\
		&\leq 
		C \cdot \bigl(\nu_{\H} \cdot \| \nablaH v_d \|_{\rL^2(\Omega)}^2 + \nu_z \cdot  \| \dz v_\rE-\dz v_{0} \|_{\rL^2(\Omega)}^2 \bigr)
		\\
		&+C \bigl(\| v_{\rE} \|_{\rL^\infty(-h,0)}^2 + \| \dz v_{\rE} \|_{\rL^\infty(-h,0)}^2 \bigr) \| v_{\rE} - v_0 \|_{\rL^2(\Omega)}^2 \\
		&+ C \cdot \int_0^t \| \partial_z v_d \|_{\rL^2(\Omega)}^2 \cdot \| \nabla \partial_z v_d \|_{\rL^2(\Omega)}^2 \cdot \| v_d \|_{\rH^1(\Omega)}^2 \\
		&+ C \cdot \int_0^t \bigl( \| \vbar_d \|_{\rH^1(\T^2)}^2 + \| \vtilde_d\|_{\rL^4(\Omega)}^4 \bigr) \cdot \| v_d \|_{\rH^1(\Omega)}^2 .
	\end{align*}
	From \autoref{lem:energy_2} and estimate \eqref{eq:bound_2} we conclude
	\begin{align*}
		\int_0^t \| \partial_z v_d \|_{\rL^2(\Omega)}^2 \cdot \| \nabla \partial_z v_d \|_{\rL^2(\Omega)}^2 
		&\leq \| \partial_z v_d \|_{\rL^\infty_t\rL^2(\Omega)}^2 \cdot  \int_0^t \| \nabla \partial_z v_d \|_{\rL^2(\Omega)}^2 
		\leq \hat{B}(t)^2, \\
		\int_0^t \| \vbar_d \|_{\rH^1(\T^2)}^2 + \| \vtilde_d\|_{\rL^4(\Omega)}^4
		&\leq \int_0^t \| v_d \|_{\rL^2(\Omega)}^2
		+
		\int_0^t \| \nablaH \vbar_d \|_{\rL^2(\T^2)}^2 + \int_0^t \| \vtilde_d\|_{\rL^4(\Omega)}^4 \\
		&\leq \| v_{\rE}-v_0 \|_{\rL^2(\Omega)}^2 +  2 \hat{B}(t),
	\end{align*} 
	and both quantities are integrable.
	Using 
	\begin{align*}
		K
		&:= C \cdot \bigl(\nu_{\H} \cdot \| \nablaH v_d \|_{\rL^2(\Omega)}^2 + \nu_z \cdot  \| \dz v_\rE-\dz v_{0} \|_{\rL^2(\Omega)}^2 \bigr) \\
		&\ +C \cdot \bigl(1+\| v_{\rE} \|_{\rL^\infty(-h,0)}^2 + \| \dz v_{\rE} \|_{\rL^\infty(-h,0)}^2 \bigr) \cdot \| v_{\rE} - v_0 \|_{\rL^2(\Omega)}^2 . 
	\end{align*}
	and applying Gronwall's inequality yields the estimate
	\begin{equation*}
		\| v_d(t) \|_{\rH^1(\Omega)}^2
		+ \int_0^t \| v_d \|_{\rH^2(\Omega)}^2
		\leq K \cdot \mathrm{e}^{\int_0^t C \cdot (1+\hat{B}(s))^2 \d s}
		=: B(t) . \qedhere 
	\end{equation*}
	%
	%
	%
	%
	%
	%
	%
	%
\end{proof}

%

In the next step we establish higher-order a-priori estimates, which proves \autoref{prop:Hk_2}.

\begin{prop}[$\rL^\infty_t\rH^k$-$\rL^2_t \rH^{k+1}$-estimates]\label{prop:estimate_2 Hk}
	Suppose $\alpha < 0$ 
	and fix $k \in \N$. 
	Let $v_d$ be a solution to \eqref{eq:difference_2}. Then there exists a bound $B_k$ continuously depending on $t, \| v_0 \|_{\rH^k(\Omega)}$ and $\| v_{\rE} \|_{\rH^{k+1}(\Omega)}$ such that 
	\begin{equation*}
		\max_{t \in (0,T)} \| v_d(t) \|_{\rH^k(\Omega)}^2
		+ \int_0^t \| v_d(s) \|_{\rH^{k+1}(\Omega)}^2 \d s \leq B_k(t) .
	\end{equation*}
\end{prop}
\begin{proof}
	We argue by induction over $k \in \N$.
	The initial case $k = 1$ is proven in \autoref{prop:estimate_2}. 
	
	For the induction step we use the maximal regularity of the hydrostatic Ekman-Stokes operator $A$ on $\rH^k(\Omega) \cap \rL^2_{\sigmabar}(\Omega)$, see \autoref{prop:mr_2}, and the embedding $\E_1^k(T) \hookrightarrow \mathrm{BUC}(0,T;\rH^{k+1}(\Omega) \cap \rL^2_{\sigmabar}(\Omega))$, and obtain
	\begin{equation*}
		\| v_d(t) \|_{\rH^{k+1}(\Omega)}^2
		+ 
		\int_0^t \| v_d(s) \|_{\rH^{k+2}(\Omega)}^2 \d s 	
		\leq 
		C \cdot \| v_0 \|_{\rH^{k+1}(\Omega)}^2 
		+ C \cdot \int_0^t \| F(v_d,v_d) \|_{\rH^k(\Omega)}^2 .
	\end{equation*} 
	Applying \autoref{lem:bilinear bounds_2} and Young's inequality we obtain
	\begin{align*}
		\| F(v_d,v_d) \|_{\rH^k(\Omega)}^2
		&\leq C \cdot \| v_d \|_{\rH^{k+2}(\Omega)} \cdot \| v_d \|_{\rH^{k+1}(\Omega)}^2 \cdot \| v_d \|_{\rH^{k}(\Omega)} \\
		&\leq \frac{1}{2} \cdot \| v_d \|_{\rH^{k+2}(\Omega)}^2
		+ C
		\cdot \| v_d \|_{\rH^{k+1}(\Omega)}^4 \cdot \| v_d \|_{\rH^{k}(\Omega)}^2 .
	\end{align*}
	Integrating in time and absorbing the term  $\int_0^t \| v_d \|_{\rH^{k+2}(\Omega)}^2$ yields
	\begin{align*}
		\| v_d(t) \|_{\rH^{k+1}(\Omega)}^2
		&+ 
		\int_0^t \| v_d(s) \|_{\rH^{k+2}(\Omega)}^2 \d s \\	
		&\leq 
		C \cdot \| v_0 \|_{\rH^{k+1}(\Omega)}^2 
		+ C \cdot \int_0^t \| v_d \|_{\rH^{k+1}(\Omega)}^2 \cdot \| v_d \|_{\rH^{k}(\Omega)}^2 \cdot \| v_d \|_{\rH^{k+1}(\Omega)}^2 .
	\end{align*} 
	We set 
	\begin{equation*}
		K_k(t) := C \cdot \| v_d(t) \|_{\rH^{k+1}(\Omega)}^2 \cdot \| v_d(t) \|_{\rH^{k}(\Omega)}^2
	\end{equation*}
	and note that by induction hypothesis we have
	\begin{equation*}
		\int_0^t K_k(s) \, \d s
		\leq C \cdot \sup_{t} \| v_d(t) \|_{\rH^{k}(\Omega)}^2
		\cdot \int_0^t  \| v_d(s) \|_{\rH^{k+1}(\Omega)}^2 \, \d s
		\leq C \cdot B_k(t)^2  
	\end{equation*}
	and $K_k \in \rL^1(0,T)$. Now Gronwall's inequality implies the desired estimate
	\begin{equation*}
		\| v_d(t) \|_{\rH^{k+1}(\Omega)}^2
		+ 
		\int_0^t \| v_d(s) \|_{\rH^{k+2}(\Omega)}^2 \d s 	
		\leq 
		C \cdot \| v_0 \|_{\rH^{k+1}(\Omega)}^2
		\cdot \mathrm{e}^{ \int_0^t C \cdot B_k(s)^2 \, \d s} =: B_{k+1}(t). \qedhere 
	\end{equation*}
\end{proof}

\begin{rem}
	In contrast to the energy estimates, \autoref{lem:energy_2}, where the right-hand side is independent of $t$ and $T$, the right-hand sides of the estimates in \autoref{prop:estimate_2 Hk} dependent on $t$ or $T$. Therefore, in \autoref{lem:energy_2} it is possible to pass to the limit $T \to + \infty$, whereas in \autoref{prop:estimate_2 Hk} it is not. In this sense the energy estimates are global estimates whereas the bounds in \autoref{prop:estimate_2 Hk} only almost global estimates.
\end{rem}

Since we now have the $\rL^2_t\rH^2$-norm under control we conclude our first global existence result.

\begin{prop}[Almost global wellposedness]\label{prop:almost global_2}
	Assume $\alpha < 0$. 
	We fix an arbitrary $T \in (0,+\infty)$. 
	Then, for every initial data $v_0 \in \rH^1(\Omega)\cap \rL^2_{\sigmabar}(\Omega)$ there exists a unique global solution $v_d \in \E_1(T) =  \rH^{1}(0,T;\rL^2_{\sigmabar}(\Omega)) \cap \rL^{2}(0,T;\rH^{2}(\Omega))$ to \eqref{eq:difference_2}.	
\end{prop}
\begin{proof}
	\autoref{prop:estimate_2} rules out the blow-up scenario from \autoref{cor:blow-up_2}. Hence, there exists a unique global solution $v_d \in \E_1(T)$.
\end{proof}

Combining almost global wellposedness with stability and the global energy estimate implies global wellposedness and convergence, which proves our main \autoref{thm:main_2}. 

\begin{thm}[Global wellposedness and convergence]\label{thm:global_2}
	Assume $\alpha < 0$. 
	Then, for every initial data $v_0 \in \rH^1(\Omega)\cap \rL^2_{\sigmabar}(\Omega)$ there exists a unique global solution $v_d \in \E_1=  \rH^{1}(\R_+;\rL^2_{\sigmabar}(\Omega)) \cap \rL^{2}(\R_+;\rH^{2}(\Omega))$ to \eqref{eq:difference_2}. Furthermore, this solution decays exponentially, so
	\begin{equation*}
		\| v_d(t) \|_{\rH^1(\Omega)} \leq C \mathrm{e}^{\omega_0 t}
		\qquad \text{ as } \enspace t \to +\infty. 
	\end{equation*} 	
	with rate $\omega_0 < 0$.
\end{thm}
\begin{proof}[Proof of \autoref{thm:main_2}.]
	From \autoref{lem:energy_2} and \autoref{prop:estimate_2} we know that 
	\begin{equation*}
		[0,+\infty) \to \R \colon t \mapsto \| v_d(t) \|_{\rH^1(\Omega)}
	\end{equation*}
	is $2$-integrable and continuous. 
	Hence $\inf_{t \geq 0} \| v_d(t) \|_{\rH^1(\Omega)} = 0$ and we may choose $t_0 \in [0,\infty)$ such that $\| v_d(t_0) \|_{\rH^1(\Omega)} <r$, where $r > 0$ is the constant from \autoref{cor:small data_2}.
	Now, we consider $\hat{v}_d(t) := v_d(t+t_0)$. Then $\hat{v}_d$ satisfy \eqref{eq:difference_2} with initial datum $\hat{v}_d(0) = v(t_0)$.
	\autoref{cor:small data_2} imply  global wellposedness of \eqref{eq:difference_2}, and hence of \eqref{eq:pe_2}, as well as
	\begin{equation*}
		\| v(t) - v_{\rE} \|_{\rH^1(\Omega)}
		= 
		\| v_d(t) \|_{\rH^1(\Omega)}
		= \| \hat{v}_d(t-t_0) \|_{\rH^1(\Omega)}
		\leq C e^{-\omega t_0} e^{\omega t}
		\leq C e^{\omega t}
	\end{equation*}  
	for $t \geq t_0$. Moreover, we know from \autoref{prop:estimate_2} that there exists a constant $C = C(t_0) > 0$ such that 
	\begin{equation*}
		\| v(t) - v_{\rE} \|_{\rH^1(\Omega)}
		= 
		\| v_d(t) \|_{\rH^1(\Omega)}
		\leq C 
	\end{equation*}
	for $0 \leq t \leq t_0$. Combining both estimates we obtain the desired decay estimate
	\begin{equation*}
		\| v(t) - v_{\rE} \|_{\rH^1(\Omega)}
		= 
		\| v_d(t) \|_{\rH^1(\Omega)}
		\leq C e^{\omega t}
	\end{equation*}  
	for all $t \geq 0$. 
	%
	%
\end{proof}

%
%


\appendix

\section{The Ekman spiral of finite depth}
\label{sec:ekman explicit_2}

This section is dedicated to providing an explicit representation of the Ekman layer. Since our depth $h >0$ is finite the representation becomes more complicated than in the infinite case. 
First, we show \autoref{rem:ekman_2}. 

\begin{lem}\label{lem:representation ekman_2}
	The Ekman spiral $v_{\rE} = (v_1,v_2)^\top$ is given by
	\begin{align*}
		v_1(z)
		&= k_1 \sin(\tfrac{z}{d}) e^{-\tfrac{z}{d}}
		+ k_2 \cos(\tfrac{z}{d}) e^{-\tfrac{z}{d}}
		+ k_3 \sin(\tfrac{z}{d}) e^{\tfrac{z}{d}}
		+ k_4 \cos(\tfrac{z}{d}) e^{\tfrac{z}{d}}, \\
		v_2(z)
		&=
		k_1 \cos (\tfrac{z}{d}) e^{-\tfrac{z}{d}}
		- k_2 \sin (\tfrac{z}{d}) e^{-\tfrac{z}{d}}
		- k_3 \cos (\tfrac{z}{d}) e^{\tfrac{z}{d}}
		+ k_4 \sin (\tfrac{z}{d}) e^{\tfrac{z}{d}},
	\end{align*}
	where $d := \sqrt{\frac{2 \nu_z}{|f|}}$ is the thickness of the Ekman layer and 
	the coefficients are given by
	\begin{equation}
		\begin{aligned}
			k_1 &:=	
			\frac{1}{2 \left(e^{\tfrac{4 h}{d}}+2 e^{\tfrac{2 h}{d}} \cos (\tfrac{2 h}{d})+1\right)}
			\cdot 
			\biggl( &&-2 e^{\tfrac{h}{d}} \bigl(e^{\frac{2 h}{d}}-1\bigr) \sin (\tfrac{h}{d})\cdot v_{g,1} \\
			& &&+2 e^{\tfrac{h}{d}} \bigl(e^{\tfrac{2 h}{d}}+1\bigr) \cos (\tfrac{h}{d}) \cdot v_{g,2}
			\\
			& &&+d \bigl(e^{\frac{2 h}{d}} \sin (\tfrac{2 h}{d})+e^{\tfrac{2 h}{d}} \cos (\tfrac{2 h}{d})+1\bigr) \cdot \tau_1 \\
			& &&-d \bigl(-e^{\tfrac{2 h}{d}} \sin (\tfrac{2 h}{d})+e^{\tfrac{2 h}{d}} \cos (\tfrac{2 h}{d})+1\bigr) \cdot \tau_2\biggr) , \\
			k_2 &:=	
			\frac{1}{2 \left(e^{\tfrac{4 h}{d}}+2 e^{\tfrac{2 h}{d}} \cos (\tfrac{2 h}{d})+1\right)}
			\cdot 
			\biggl(
			&&  -2 e^{\tfrac{h}{d}} \bigl(e^{\tfrac{2 h}{d}}+1\bigr) \cos (\tfrac{h}{d}) \cdot v_{g,1} \\
			& &&-2 e^{\tfrac{h}{d}} \bigl(e^{\tfrac{2 h}{d}}-1\bigr) \sin (\tfrac{h}{d}) \cdot v_{g,2}
			\\
			& &&+d \bigl(-e^{\tfrac{2 h}{d}} \sin (\tfrac{2 h}{d})+e^{\tfrac{2 h}{d}} \cos (\tfrac{2 h}{d})+1\bigr) \cdot \tau_1 \\
			& &&+d \bigl(e^{\tfrac{2 h}{d}} \sin (\tfrac{2 h}{d})+e^{\tfrac{2 h}{d}} \cos (\tfrac{2 h}{d})+1\bigr) \cdot \tau_2 \biggr), \\
				k_3 &:=	
			\frac{1}{2 \left(e^{\tfrac{4 h}{d}}+2 e^{\tfrac{2 h}{d}} \cos (\tfrac{2 h}{d})+1\right)}
			\cdot e^{\tfrac{h}{d}}
			\biggl( \hspace{-1 em}
			&&\ 2 \bigl(e^{\tfrac{2 h}{d}}-1\bigr) \sin (\tfrac{h}{d}) \cdot v_{g,1} \\
			& &&-2 \bigl(e^{\tfrac{2 h}{d}}+1\bigr) \cos (\tfrac{h}{r}) \cdot v_{g,2}
			\\
			& &&+d e^{\tfrac{h}{d}} \bigl(e^{\tfrac{2 h}{d}}-\sin (\tfrac{2 h}{d})+\cos (\tfrac{2 h}{d})\bigr) \cdot \tau_1 \\
			& &&-d e^{\tfrac{h}{d}} \bigl(e^{\tfrac{2 h}{d}}+\sin (\tfrac{2 h}{d})+\cos (\tfrac{2 h}{d})
			\bigr) \cdot \tau_2 \biggr), 
		\end{aligned}			\label{eq:coefficients_2}
		\end{equation}
		\begin{equation*}
			\begin{aligned}
			k_4 &:=	
			\frac{1}{2 \left(e^{\tfrac{4 h}{d}}+2 e^{\tfrac{2 h}{d}} \cos (\tfrac{2 h}{d})+1\right)}
			\cdot e^{\tfrac{h}{d}}
			\biggl( \hspace{-1 em}
			&&\ 2 \bigl(e^{\tfrac{2 h}{d}}+1\bigr) \cos (\tfrac{h}{d}) \cdot v_{g,1} \\
			& &&+ 2 \bigl(e^{\frac{2 h}{d}}-1\bigr) \sin (\tfrac{h}{d}) \cdot v_{g,2}
			\\
			& &&+d e^{\tfrac{h}{d}} \bigl(e^{\tfrac{2 h}{d}}+\sin (\tfrac{2 h}{d})+\cos (\tfrac{2 h}{d})\bigr) \cdot \tau_1 \\
			&
			&&+d e^{\tfrac{h}{d}} \bigl(e^{\tfrac{2 h}{d}}-\sin (\tfrac{2 h}{d})+\cos (\tfrac{2 h}{d})
			\bigr) \cdot \tau_2 \biggr).  
		\end{aligned}
	\end{equation*}
\end{lem}
\begin{proof}
	Direct calculations.
\end{proof}

We conclude the following formulas for the derivative of the Ekman spiral. 

\begin{lem}\label{lem:derative ekman_2}
	For $v_{\rE} = (v_1,v_2)^\top$ we have
	\begin{align*}
		\partial_z v_1(z) &= 
		\frac{1}{d} \cdot \biggl(
		k_1 \cdot \bigl( \cos(\tfrac{z}{d}) -\sin(\tfrac{z}{d}) \bigr) e^{-\tfrac{z}{d}}
		+
		k_2 \cdot \bigl( -\sin(\tfrac{z}{d}) -\cos(\tfrac{z}{d}) \bigr) e^{-\tfrac{z}{d}} \\
		&\qquad +
		k_3 \cdot \bigl( \sin(\tfrac{z}{d}) + \cos(\tfrac{z}{d}) \bigr) e^{\tfrac{z}{d}}
		+
		k_4 \cdot \bigl( -\sin(\tfrac{z}{d}) +\cos(\tfrac{z}{d}) \bigr) e^{\tfrac{z}{d}}
		\biggr),
		\\
		\partial_z v_2(z) &= 
		\frac{1}{d} \cdot \biggl(
		k_1 \cdot \bigl( -\sin(\tfrac{z}{d}) -\cos(\tfrac{z}{d}) \bigr) e^{-\tfrac{z}{d}}
		+
		k_2 \cdot \bigl( \sin(\tfrac{z}{d}) -\cos(\tfrac{z}{d}) \bigr) e^{-\tfrac{z}{d}} \\
		&\qquad +
		k_3 \cdot \bigl( \sin(\tfrac{z}{d}) -\cos(\tfrac{z}{d}) \bigr) e^{\tfrac{z}{d}}
		+
		k_4 \cdot \bigl(\sin(\tfrac{z}{d}) +\cos(\tfrac{z}{d}) \bigr) e^{\tfrac{z}{d}}
		\biggr),
	\end{align*}
	where $d := \sqrt{\frac{\nu_z}{2|f|}}$ is the thickness of the Ekman layer and
	the coefficients $k_1,\dots,k_4  \in \R$ are given by \eqref{eq:coefficients_2}.
\end{lem}
\begin{proof}
	Direct calculations.
\end{proof}

\subsection*{Acknowledgement}

The author wants to thank Felix Brandt and Matthias Hieber for fruitful discussions and helpful comments about the topic. Furthermore, the author is deeply grateful to the anonymous referee for the insightful suggestion to replace a technical smallness condition of the parameters in the earlier version with an operator-theoretic one, leading to the complete classification in the main result. Special thanks go to Felix Brandt for his help with the résumé. 

\subsection*{Funding}

The author would like to thank DFG for support through the project ''Globale Existenz und Singularitäten geophysikalischer Flüsse'', project no. 538212014.

\bigskip 

\end{document}